\documentclass[12pt,a4paper]{article}
\usepackage{amsmath}
\usepackage{amssymb}
\usepackage{amsthm}
\usepackage{amsfonts}
\usepackage{latexsym}

\theoremstyle{plain}
\newtheorem{thm}{Theorem}
\newtheorem{cor}{Corollary}
\newtheorem{lem}{Lemma}
\newtheorem{prop}{Proposition}
\theoremstyle{definition}
\newtheorem{defn}{Definition}

\theoremstyle{remark}
\newtheorem{claim}{Claim}[section]

\title{The Szeg\"{o} kernel of a symplectic quotient}
\author{Roberto Paoletti\footnote{\noindent{\bf Address.}
Dipartimento di Matematica e Applicazioni, Universit\`a degli
Studi di Milano Bicocca, Via Bicocca degli Arcimboldi 8, 20126
Milano, Italy; {\bf e-mail}: roberto.paoletti@unimib.it }}
\date{}

\begin{document}
\maketitle
\section{Introduction.}

The object of this paper is the relation between the Szeg\"{o}
kernel of an ample line bundle on a complex projective manifold,
$M$, and the Szeg\"{o} kernel of the induced polarization on the
quotient of $M$ by the holomorphic action of a compact Lie group,
$G$.

Let $M$ be an $n$-dimensional complex projective manifold and $L$
an ample line bundle on it. Suppose a connected compact Lie group
$G$ acts on $M$ as a group of holomorphic automorphisms, and that
the action linearizes to $L$. Without loss of generality, we may
then choose a $G$-invariant K\"{a}hler form $\Omega$ on $M$
representing $c_1(L)$. We may also assume given a $G$-invariant
Hermitian metric $h$ on $L$ such that the unique Hermitian
connection on $L$ compatible with the holomorphic structure has
normalized curvature $-2\pi i\Omega$.

Let $\Phi :M\rightarrow \frak{g}^*$, where $\frak{g}$ is the Lie
algebra of $G$, be the moment map for the action (this is
essentially equivalent to the assignment of a linearization).
Suppose that $0\in \frak{g}^*$ lies in the image of $\Phi$, that
it is a regular value of $\Phi$, and furthermore that $G$ acts
freely on $\Phi ^{-1}(0)$. In this situation, the Hilbert-Mumford
quotient for the action of the complexification $\tilde G$ on $M$
is nonsingular, and may be naturally identified with the
symplectic reduction $M_0=\Phi ^{-1}(0)/G$. The latter,
furthermore, inherits a naturally induced polarization $L_0$, with
an Hermitian metric $h_0$ and a K\"{a}hler form $\Omega _0$. The
quotient structures $L_0,h_0,\Omega _0$ are induced simply by
descending the restrictions of $L,h,\Omega$ to $\Phi ^{-1}(0)$
down to $M_0$ \cite{gs-gq}.

Let $H^0(M,L^{\otimes k})$ and $H^0(M_0,L^{\otimes k}_0)$ denote
the spaces of holomorphic sections of powers of $L$ on $M$ and of
powers of $L_0$ on $M_0$, respectively. Then $G$ acts linearly on
$H^0(M,L^{\otimes k})$, and in the given hypothesis by the theory
of \cite{gs-gq} there is for every integer $k\ge 0$ a natural
isomorphism $H^0(M,L^{\otimes k})^G\cong H^0(M_0,L^{\otimes
k}_0)$. Here $H^0(M,L^{\otimes k})^G\subseteq H^0(M,L^{\otimes
k})$ is the subspace of $G$-invariant holomorphic sections.

The given choices equip the vector spaces  $H^0(M,L^{\otimes k})$
and $H^0(M_0,L^{\otimes k}_0)$ with the following unitary
structures $(\,,\,)$ and $(\,,\,)_0$, respectively. If $\sigma,
\tau \in H^0(M,L^{\otimes k})$, then
$$(\sigma,\tau)\,=:\,\int _Mh_p^{\otimes k}(\sigma (p),\tau(p))
\mathrm{vol}_M(p),$$ where $\mathrm{vol}_M=:\Omega ^{\wedge n}$ is
the volume form on $M$ associated to the K\"{a}hler form $\Omega$,
and $h^{\otimes k}$ denotes the induced Hermitian metric on
$L^{\otimes k}$. Similarly, if $\sigma_0, \tau _0\in
H^0(M_0,L^{\otimes k}_0)$, then
$$(\sigma_0,\tau_0)_0\,=:\,\int _{M_0}h_{0,p}^{\otimes k}(\sigma _0(p),\tau_0(p))
\mathrm{vol}_{M_0}(p),$$ where now $\mathrm{vol}_{M_0}=:\Omega_0
^{\wedge (n-g)}$ is the volume form on $M_0$ associated to the
K\"{a}hler form $\Omega_0$.

It is then natural to ask the following question: what are the
asymptotic metric properties, with respect to these unitary
structures, of the isomorphisms $r_k:H^0(M,L^{\otimes
k})^G\rightarrow H^0(M_0,L^{\otimes k}_0)$ as $k\rightarrow
+\infty$? To simplify the exposition, we shall mostly leave the
isomorphisms $r_k$ implicit. From an analytic viewpoint, this
question may be phrased as follows: how does the Szeg\"{o} kernel
of the pair $(M,L)$ relate to the Szeg\"{o} kernel of the quotient
pair $(M_0,L_0)$?

A natural measure of the difference between the two unitary
structures is offered by the comparison of orthonormal basis. In
this direction, recall the following basic fact from the theory of
algebro-geometric Szeg\"{o} kernels \cite{z}: given for every
integer $k\ge 0$ an orthonormal basis $\left \{t_j^{(k)}\right \}$
of $H^0(M_0,L_0^{\otimes k})$, we have an asymptotic expansion,
uniform in $p_0\in M_0$,
\begin{equation}\label{eqn:basicexpansion}
\sum _j\left \|t_j^{(k)}(p_0)\right \|^2\,\sim\,k^{n-g}+\sum
_{l\ge 1}b_l(p_0)k^{n-g-l},\end{equation} where $n$ is the complex
dimension of $M$, $g$ the real dimension of $G$, and the $a_l$'s
are smooth functions on $M_0$; the right hand side of
(\ref{eqn:basicexpansion}) is of course independent of the
particular choice of the orthonormal basis $\left
\{t_j^{(k)}\right \}$. In order to measure how far an orthonormal
basis of $H^0(M,L^{\otimes k})^G$ is from being an orthonormal
basis of $H^0(M_0,L_0^{\otimes k})$, we may then look for a
similar asymptotic expansion for orthonormal basis of
$H^0(M,L^{\otimes k})^G$, along $\Phi ^{-1}(0)\subseteq M$.

The following is a special case of Theorem 1 of \cite{p}: Given
for every integer $k\ge 0$ an orthonormal basis $\left
\{s_j^{(G,k)}\right \}$ of $H^0(M,L^{\otimes k})^G$, we have an
asymptotic expansion, uniform in $p\in \Phi ^{-1}(0)$:
\begin{equation}\label{eqn:basicexpansionG}
\sum _j\left \|s_j^{(G,k)}(p)\right
\|^2\,\sim\,a_0(p)\,k^{n-g/2}+\sum _{l\ge
1}a_l(p)k^{n-g/2-l},\end{equation} where the $a_l$'s are smooth
$G$-invariant functions on $\Phi ^{-1}(0)$, and $a_0(p)>0$ for
every $p\in \Phi ^{-1}(0)$. The same expression is rapidly
decaying when $p\not\in \Phi ^{-1}(0)$. We are thus led to ask
whether the $r_k$'s are asymptotically conformally isometric, by a
conformal factor involving an appropriate power of $k$.

The following Theorem shows that this is not the case, unless the
effective potential $V_{\mathrm{eff}}$ of the action is a
constant. Recall that $V_{\mathrm{eff}}$ is the smooth function on
$M_0$ whose value at $p_0\in M_0$ is the volume of the fibre $\Phi
^{-1}(p_0)\subseteq \Phi^{-1}(0)$ in the restricted metric. This
is an important attribute of the action, playing a crucial role in
the study of the K\"{a}hler structure of the quotient
\cite{burnsg}. It may be viewed as the $G$-invariant function on
$\Phi ^{-1}(0)$ associating to each $p\in \Phi ^{-1}(0)$ the
volume of the fibre $G\cdot p$.

\begin{thm}\label{thm:mainexpgequiv}
Suppose as above that $0\in \frak{g}^*$ is a regular value of the
moment map, and that $G$ acts freely on $\Phi ^{-1}(0)$. Let $n$
be the complex dimension of $M$ and $g$ be the real dimension of
$G$. For every integer $k\ge 0$ let $\left \{s_j^{(G,k)}\right \}$
be an orthonormal basis of the space $H^0(M,L^{\otimes k})^G$ of
$G$-invariant holomorphic sections of $L^{\otimes k}$. Then there
is an asymptotic expansion, uniform in $p\in \Phi^{-1}(0)$,
\begin{equation}\label{eqn:basicexpansionGnew}
\sum _j\left \|s_j^{(G,k)}(p)\right
\|^2\,\sim\,\frac{1}{V_{\mathrm{eff}}(p)}\,k^{n-g/2}+\sum _{l\ge
1}a_l(p)k^{n-g/2-l}.\end{equation} \end{thm}

Let now the weights $\{\omega\}$ label the finite dimensional
irreducible representations of the connected compact Lie group
$G$; let $V_\omega$ be the representation associated to the weight
$\omega$. Given the linear action of $G$ on $H^0(M,L^{\otimes
k})$, there is a $G$-equivariant direct sum decomposition
$$H^0(M,L^{\otimes k})\, \cong\, \bigoplus _{\omega}H^0(M,L^{\otimes
k})_\omega$$ of $H^0(M,L^{\otimes k})$ over the finite dimensional
irreducible representations of $G$; for every weight $\omega$ the
subspace $H^0(M,L^{\otimes k})_\omega$ is $G$-equivariantly
isomorphic to a direct sum of copies of $V_\omega$. Pairing
Theorem \ref{thm:mainexpgequiv} above with Theorem 1 of \cite{p},
we immediately get:

\begin{cor}\label{cor:expgeneralrepres}
In the hypothesis of Theorem \ref{thm:mainexpgequiv}, let $\omega$
be the weight corresponding to a finite dimensional irreducible
representation of $G$. For every integer $k\ge 0$ let $\left
\{s_j^{(\omega,k)}\right \}$  be an orthonormal basis of the
subspace $H^0(M,L^{\otimes k})_\omega\subseteq H^0(M,L^{\otimes
k})$ . Then there is an asymptotic expansion, uniform in $p\in
\Phi^{-1}(0)$,
\begin{equation}\label{eqn:basicexpansionGomega}
\sum _j\left \|s_j^{(\omega,k)}(p)\right \|^2\,\sim\,\frac{\dim
(V_\omega)^2 }{V_{\mathrm{eff}}(p)}\,k^{n-g/2}+\sum _{l\ge
1}a_{l,\omega}(p)k^{n-g/2-l}.\end{equation} \end{cor}

The asymptotic expansion in Theorem \ref{thm:mainexpgequiv} is a
straightforward consequence of the following analytic core result,
which we state here rather loosely and that will be spelled out
more precisely in the course of the paper:

\bigskip
\noindent Key Result: \textit{The Szeg\"{o} kernel of the triple
$(M,L,h)$ descends, in an appropriate sense, to an elliptic
Toeplitz operator on the circle bundle of the quotient triple
$(M_0,L_0,h_0)$, and after a suitable renormalization the symbol
of this Toeplitz operator is (proportional to) the effective
volume of the action.}

\bigskip

The techniques in this paper are based in a general sense on the
microlocal theory of the Szeg\"{o} kernel \cite{bs}, and more
specifically on its formulation based on Fourier-Hermite
distributions developed in \cite{bg}; in this framework one has a
good control of the functorial behaviour of Toeplitz operators and
their symbols under geometric operations like restriction and
push-forward. The general strategy used here was also inspired by
the approach of Shiffman, Tate and Zelditch to Szeg\"{o} kernels
on toric varieties \cite{stz}, especially by their philosophy of
restricting the Szeg\"{o} kernel of an ambient projective space to
a projective submanifold. In fact, here we shall first restrict
the Szeg\"{o} kernel of $(M,L,h)$ to the locus where the moment
map vanishes, and then push it forward by the $G$-action.

\section{Metalinear and metaplectic preliminaries.}
We shall make use of the notions of metalinear manifold and
half-form \cite{gs-ga}. If $Z$ is a manifold, we shall denote by
$\left | \bigwedge \right |(Z)$ the line bundle of densities on
$Z$, and by $\left | \Omega \right |(Z)$ its space of smooth
global sections, which we simply call densities on $Z$. If $Z$ is
a metalinear manifold, we shall denote by $\bigwedge ^{1/2}(Z)$
the line bundle of metalinear forms on $Z$, and by $\Omega
^{1/2}(Z)$ its space of smooth global sections, which we simply
call metalinear forms on $Z$. An orientation on a manifold induces
a metalinear structure on it, and a volume form induces a nowhere
vanishing half-form, its \textsl{square root}. If $\mu, \, \nu\in
\Omega ^{1/2}(Z)$, the product $\mu \cdot \overline{\nu}$ is a
density on $X$.  The space of square integrable half-forms on an
metalinear manifold has a natural Hilbert structure, given by
$$<\,\mu,\,\nu\,>\,=\, \int _Z\, \mu \cdot
\overline{\nu}\,\,\,\,\,\,\,(\mu,\,\nu \,\in \,\Omega
^{1/2}(Z)).$$

There is a more general notion of metalinear vector bundle, a
metalinear manifold being a manifold with a metalinear tangent
bundle. To a metalinear structure on a rank-$r$ vector bundle $E$
there is associated a line bundle $\bigwedge ^{1/2}(E)$ satisfying
$$\bigwedge ^{1/2}(E)\otimes \bigwedge ^{1/2}(E)\cong \bigwedge
^r(E^*),$$ where the latter is the top exterior power of the dual
$E^*$. The space of smooth global sections of $\bigwedge
^{1/2}(E)$ will be denoted by $\Omega ^{1/2}_{E}$. Thus, if $Z$ is
a metalinear manifold then $\bigwedge ^{1/2}(Z)=\bigwedge
^{1/2}(TZ)$ and $\Omega ^{1/2}(Z)=\Omega ^{1/2}_{TZ}$. If $E$ is
oriented, it is metalinear, and if $0\rightarrow A\rightarrow
E\rightarrow B\rightarrow 0$ is a short exact sequence of vector
bundles on a manifold, then any of the three vector bundles is
metalinear if the other two are, and $\bigwedge ^{1/2}(E)\cong
\bigwedge ^{1/2}(A)\otimes \bigwedge ^{1/2}(B)$. A nowhere
vanishing smooth section of $\bigwedge ^{r}(E^*)$ induces a
metalinear structure and a nowhere vanishing section of $\bigwedge
^{1/2}_{E}$, $\sqrt \sigma \in \Omega ^{1/2}_{E}$.

In the following, we shall be making various natural choices of a
nowhere vanishing section of $\bigwedge ^{r}(E^*)$ for some vector
bundles in our picture ($r$ being the rank of $E$); we shall
generally denote this choice by $\mathrm{vol}_E$ and the
associated nowhere vanishing section of $\bigwedge ^{1/2}(E)$ by
$\mathrm{vol}^{(1/2)}_E=\sqrt{\mathrm{vol}_E}$.

Let $P$ and $Q$ be metalinear manifolds. A morphism of metalinear
manifolds from $P$ to $Q$ is the assignment of a smooth map
$f:P\rightarrow Q$ and a morphism of vector bundles $\tilde f
:f^*\left ( \bigwedge ^{1/2}(Z)\right )\rightarrow \bigwedge
^{1/2}(P)$. Let $f:P\rightarrow Q$ be any smooth map. Let
$$\Gamma _f=\left \{\,\big (\,
(p,-(d_pf)^t(\eta)),\,(f(p),\eta)\,)\,:\, p\in P,\, \eta \in
T_{f(p)}^*Q\right \}$$ be the conormal bundle to the graph of $f$.
Then $\Gamma _f$ is clearly a Lagrangian submanifold of
$T^*P\times T^*Q$ and a metalinear manifold. Furthermore, giving
$f$ the structure of a morphism of metalinear manifolds is
equivalent to assigning an appropriate half-form on $\Gamma _f$
\cite{gs-ga}.

Let us suppose, in particular, that $Q$ is an oriented Riemannian
manifold. We shall denote by $\eta _Q$ the half-form taking value
one on oriented orthonormal frames of $TQ$, and call it the
\textit{canonical} half-form of $Q$. If $P$ is also an oriented
Riemannian manifold, we can then make any smooth map $f
:P\rightarrow Q$ into a morphism of metalinear manifolds by
setting $\tilde f(\eta _Q)=\eta _P$. Let $\Gamma _f\subseteq
T^*(P\times Q)\setminus \{0\}$ be the conormal bundle to the graph
of $f$, with projections $q_1:\Gamma _f \rightarrow P$ and
$q_2:\Gamma _f \rightarrow Q$. Given the exact sequence
$$0\longrightarrow q_2^*(T^*Q)\rightarrow T(\Gamma
_f)\longrightarrow q_1^*(TP)\longrightarrow 0,$$ the half-form on
$\Gamma_f$ associated to $\tilde f$ is
\begin{equation}\label{eqn:halfformmeta}
 \mathrm{vol}_{\Gamma _f}^{1/2}=:q_2^*(\eta _Q^{-1})\otimes
q_1^*(\eta _P).
\end{equation}

If a morphism of metalinear manifolds is given, passing to global
sections we can define the pull-back of a metalinear half-form on
$Q$ to a metalinear half-form on $P$; we shall simply denote this
by $f^*:\Omega ^{1/2}(Q)\rightarrow \Omega ^{1/2}(P)$. If, in
addition, $f$ is a proper submersion, we can also define a
push-forward operation $f_*:\Omega ^{1/2}(P)\rightarrow \Omega
^{1/2}(Q)$. The microlocal theory for these operations on Fourier
and Fourier-Hermite generalized half-forms has been developped in
\cite{gs-ga} and \cite{bg}.

The proof of the following Lemma is left to the reader:

\begin{lem}
Let $P$, $Q$ be oriented Riemannian manifolds, and let
$f:P\rightarrow Q$ be a proper submersion. Let us make $f$ into a
morphism of metalinear manifolds by setting $\tilde f^*(\eta
_Q)=:\eta _P$. For $q\in Q$, let $\upsilon (q)$ denote the volume
of the fibre $f^{-1}(q)$ in the induced metric structure. Then the
associated push-forward operation $f_*:\Omega ^{1/2}(P)\rightarrow
\Omega ^{1/2}(Q)$ satisfies
$$f_*(\eta _P)=\upsilon \cdot \eta _Q.$$\label{lem:propersubmersion}
\end{lem}

We shall now review some basic facts and notation from metaplectic
geometry. Let $\mathrm{Sp}(\ell)$ be the group of $2\ell\times
2\ell$ (real) symplectic matrices, and for $1\le k\le \ell$
define:
$$
\mathrm{Sp}(k,\ell)\,=\,\left \{A=[a_{ij}]\in \mathrm{Sp}(\ell)\,
:\, a_{ij}=0\, \text{ if }\, 1\le j\le k,\, i>k\right \}.$$ There
is an obvious (surjective) morphism of Lie groups
\begin{equation}\label{eqn:liegrmorph}
\mathrm{Sp}(k,\ell)\, \stackrel{\varpi}{\longrightarrow} \,
\mathrm{GL}(k)\times \mathrm{Sp}(\ell-k).
\end{equation}
The following lemma is left to the reader:
\begin{lem}
We have:
\begin{equation*}
\begin{array}{lll}
   \ker(\varpi) & =& \left \{\begin{pmatrix}
  I_{k} & B & A_{11} & A_{12} \\
  0 & I_{\ell-k} & A_{12}^t & A_{22} \\
  0 & 0 & I_{k} & 0 \\
  0 & 0 & -B^t & I_{\ell-k}
\end{pmatrix}:B\in M_{k,\ell-k}(\mathbb{R}),\, A_{12}\in M_{k,\ell-k}(\mathbb{R}),\right.\\
     & &\left .\text{ and }A_{11}\in
M_k(\mathbb{R}),\,A_{22}\in M_{\ell -k}(\mathbb{R})\text{\textrm{
are both symmetric }}\right \}.
  \end{array}\end{equation*}
\end{lem}
Let $T=:\ker(\varpi) \subseteq \mathrm{Sp}(k,\ell)$. Then $T$ is a
contractible, hence simply connected, normal subgroup of
$\mathrm{Sp}(k,\ell)$.

Similarly, let us define
\begin{equation*}
\begin{array}{lll}
\mathrm{Sp}(k,\ell)_+&=&\left \{\begin{pmatrix}
  a_{11} & \cdots & a_{1k} & a_{1\,k+1} & \cdots & a_{1\ell} \\
  \vdots & \cdots & \vdots & \vdots & \cdots & \vdots \\
  a_{k1} & \cdots & a_{kk} & a_{k\,k+1} & \cdots & a_{k\ell} \\
  0 & \cdots & 0 & a_{k+1\,k+1} & \cdots & a_{k+1\,\ell} \\
  \vdots & \vdots & \vdots & \vdots & \vdots & \vdots\\
0 & \cdots & 0 & a_{\ell\,k+1} & \cdots & a_{\ell\,\ell} \\
\end{pmatrix}\in \mathrm{Sp}(k,\ell) \right .\\
 & &:\left. \det \begin{pmatrix}
  a_{11} & \cdots & a_{1k} \\
  \vdots & \cdots & \vdots  \\
  a_{k1} & \cdots & a_{kk} \\
\end{pmatrix}>0\right \}.\end{array}\end{equation*}
By restriction of $\varpi$, we obtain a surjective morphism of Lie
groups, with kernel $T\subseteq \mathrm{Sp}(k,\ell)_+$:
\begin{equation}\label{eqn:liegrmorph+}
\mathrm{Sp}_+(k,\ell)\, \stackrel{\varpi_+}{\longrightarrow} \,
\mathrm{GL}_+(k)\times \mathrm{Sp}(\ell-k).
\end{equation}
Therefore, $\mathrm{Sp}_+(k,\ell)$ is a connected subgroup of
$\mathrm{Sp}(\ell)$.

Let now $s:\mathrm{Mp}(\ell)\rightarrow \mathrm{Sp}(\ell)$ be the
metaplectic double cover. Now $\mathrm{GL}(\ell)$ sits naturally
in $\mathrm{Sp}(\ell)$ in the standard manner, and \begin{lem}
$s^{-1}\big (\mathrm{GL}(\ell)\big )\cong \mathrm{ML}(\ell)$.
\end{lem}
Let us next define $\mathrm{Mp}(k,\ell)=:s^{-1}\left
(\mathrm{Sp}(k,\ell)\right )$, and
$\mathrm{Mp}_+(k,\ell)=:s^{-1}\left (\mathrm{Sp}_+(k,\ell)\right
)$.
\begin{lem} $\mathrm{Mp}_+(k,\ell)$ is a connected Lie group.
\end{lem}

\noindent \textit{Proof.} Let us consider the injective morphism
of Lie groups $\psi :\mathrm{Sp}(\ell-k)\rightarrow
\mathrm{Sp}(\ell)$ given by
\begin{equation}
\begin{pmatrix}
  A_{1} & A_{2} \\
  A_{3} & A_{4}
\end{pmatrix}
\mapsto \begin{pmatrix}
  I_{k} & 0 & 0 & 0 \\
  0 & A_{1} & 0 & A_{2} \\
  0 & 0 & I_k & 0 \\
  0 & A_{3} & 0 & A_{4}
\end{pmatrix}.
\end{equation}
Clearly, $\psi \left (\mathrm{Sp}(\ell -k)\right )\subseteq
\mathrm{Sp}(k,\ell)_+$, and $\psi$ induces an isomorphism of
homotopy groups $\psi _*:\pi _1\big (\mathrm{Sp}(\ell -k)\big
)\cong \pi _1\big (\mathrm{Sp}(\ell)\big )\cong \mathbb{Z}$. Thus,
$$s^{-1}\big [\psi \big (\mathrm{Sp}(\ell -k)\big )\big ]\cong
\mathrm{Mp}(\ell -k)$$ is the unique connected 2:1 cover of
$\mathrm{Sp}(\ell-k)$. Hence, on the one hand, $\mathrm{Mp}(\ell
-k)$ is contained in the identity component of
$\mathrm{Mp}(k,\ell)_+$, and on the other if there was another
connected component of $\mathrm{Mp}(k,\ell)_+$ then points in
$\psi \left (\mathrm{Sp}(\ell -k)\right )$ would have more than
two inverse images in $\mathrm{Mp}(\ell)$, absurd.

\bigskip

Let us now consider the surjective homomorphism $\gamma
:\mathrm{Sp}(k,\ell)_+\rightarrow \mathrm{Sp}(\ell-k)$ given by
the composition of $\varpi$ with projection onto the second
factor. The following Lemma is left to the reader:

\begin{lem}
We have:
\begin{equation*}
\begin{array}{lll}
   \ker(\gamma) & =& \left \{\begin{pmatrix}
  A & B & A_{11} & A_{12} \\
  0 & I_{\ell-k} & A_{12}^t & A_{22} \\
  0 & 0 & (A^{-1})^t & 0 \\
  0 & 0 & -(A^{-1}B)^t & I_{\ell-k}
\end{pmatrix}:\right. \\
& & \left.A\in \mathrm{GL}_+(k),\, B\in M_{k,\ell-k}(\mathbb{R}),\, A_{12}\in M_{k,\ell-k}(\mathbb{R}),\right.\\
     & &\left .\text{\textrm{and }}A_{11}\in
M_k(\mathbb{R}),\,A_{22}\in M_{\ell -k}(\mathbb{R})\text{\textrm{
are both symmetric }}\right \}.
  \end{array}\end{equation*}
\end{lem}
Thus, $Z=:\ker (\gamma)$ retracts to $\mathrm{GL}_+(k)\subseteq
\mathrm{Sp}_+(k,\ell)$. Now $s^{-1}\big (\mathrm{GL}(\ell)\big
)\subseteq \mathrm{Mp}(\ell)$ is simply the metalinear group
$\mathrm{ML}(\ell)$, and consists of four connected components.
Therefore, the inverse image $s^{-1} (Z)\subseteq
\mathrm{Mp}_+(k,\ell)$ is the union of two connected components,
each mapping isomorphically onto $Z$. We shall also denote by $Z$
the identity component of $s^{-1}(Z)$.

\begin{lem}\label{lem:metapl}
$\mathrm{Mp}_+(k,\ell) /Z\cong \mathrm{Mp}(\ell-k).$
\end{lem}

\noindent
 \textit{Proof.} The commutative diagram
$$
  \begin{array}{ccc}
    \mathrm{Mp}_+(k,\ell) & \longrightarrow & \mathrm{Mp}_+(k,\ell) /Z\\
     & & \\
\downarrow & & \downarrow \\
     & & \\
   \mathrm{Sp}_+(k,\ell) & \longrightarrow & \mathrm{Sp}(\ell-k) =
\mathrm{Sp}_+(k,\ell)/Z
  \end{array}
$$
exhibits $\mathrm{Mp}_+(k,\ell) /Z$ as a connected double cover of
$\mathrm{Sp}(\ell-k)$.

\bigskip

$\mathrm{GL}(k)$ and $\mathrm{Sp}(\ell-k)$ sit naturally in
$\mathrm{Sp}(k,\ell)$, namely as the commuting Lie subgroups of
all matrices of the form
$$
\begin{pmatrix}
  A & 0 & 0 & 0 \\
  0 & I_{\ell-k} & 0 & 0 \\
  0 & 0 & (A^t)^{-1} & 0 \\
  0 & 0 & 0 & I_{\ell-k}
\end{pmatrix}
\,\,\,\,\,\,\,(A\in \mathrm{GL}(k))$$ and
$$
\begin{pmatrix}
  I_k& 0 & 0 & 0 \\
  0 & R & 0 & S \\
  0 & 0 & I_k & 0 \\
  0 & T & 0 & U
\end{pmatrix}\,\,\,\,\,\,\,\,\left (\begin{pmatrix}
 R & S \\
  T & U
\end{pmatrix}\in \mathrm{Sp}(\ell-k)\right ),
$$
respectively. We have therefore an injective morphism of Lie
groups $$o=(\zeta,\chi):\mathrm{GL}(k)\times
\mathrm{Sp}(\ell-k)\rightarrow \mathrm{Sp}(k,\ell)\subseteq
\mathrm{Sp}(\ell).$$ The inverse images $s^{-1}\big
(\mathrm{GL}(k)\big )$ and $s^{-1}\big (\mathrm{Sp}(\ell-k)\big )$
are isomorphic, respectively, to $\mathrm{ML}(k)$ and
$\mathrm{Mp}(\ell-k)$. Thus, there is a commutative diagram of Lie
group homomorphisms
$$
\begin{array}{ccccc}
    \mathrm{ML}(k) & \stackrel{\tilde \varsigma}{\longrightarrow}&\mathrm{Mp}(k,\ell)& \stackrel{\tilde \xi}{\longleftarrow}&\mathrm{Mp}(\ell-k)\\
 & & & & \\
\downarrow & &\downarrow & & \downarrow\\
 & & & & \\
 \mathrm{GL}(k) & \stackrel{\varsigma}{\longrightarrow}&\mathrm{Sp}(k,\ell)&
\stackrel{\xi}{\longleftarrow}&\mathrm{Sp}(\ell-k).
  \end{array}$$
Taking products, we have a smooth map
$$\tilde o=(\tilde \xi,\tilde \chi):\mathrm{ML}(k)\times
\mathrm{Mp}(\ell-k)\longrightarrow \mathrm{Mp}(k,\ell)\subseteq
\mathrm{Mp}(\ell),\,\,\,\,(g,h)\mapsto \tilde \xi (g)\cdot \tilde
\chi (h).$$

\begin{lem} $\tilde o$ is a Lie group homomorphism.
\end{lem}

\noindent \textit{Proof.} Given that $\tilde \xi$ and $\tilde
\chi$ are morphisms, it suffices to show that
$$\tilde \xi (g)\cdot \tilde \chi (h)\cdot \tilde \xi
(g)^{-1}\cdot \tilde \chi (h)^{-1}\,=\,e,\,\,\,\,\forall\, (g,h)
\in \mathrm{ML}(k)\times \mathrm{Mp}(\ell-k).$$ Given that $o$ is
a morphism, these commutators all lie in $\ker (s)$, which is
finite of order $2$. Given that each of $\tilde \xi$ and $\tilde
\chi$ is a morphism, they all equal $e$ if at least one of $g$ and
$h$ is the identity. By connectedness of $\mathrm{Mp}(\ell-k)$,
the claim follows.

\bigskip
The composition
$$\varpi \circ \left .s\right
|_{\mathrm{Mp}(k,\ell)}\circ \tilde o:\mathrm{ML}(k)\times
\mathrm{Mp}(\ell-k)\longrightarrow \mathrm{GL}(k)\times
\mathrm{Sp}(\ell-k)$$ is clearly the product of the double covers
$\mathrm{ML}(k)\rightarrow \mathrm{SL}(k)$ and $\mathrm{Mp}(\ell
-k)\rightarrow \mathrm{Sp}(\ell-k)$, and is thus a 4:1 covering.

\bigskip

Being contractible, $T$ lifts isomorphically to a subgroup of
$\mathrm{Mp}(\ell)$, that we still denote  $T\lhd
\mathrm{Mp}(k,\ell)$. Let us define the quotient group
$\mathrm{Mp}(k|\ell)=\mathrm{Mp}(k,\ell)/T$. Let
$\theta:\mathrm{Mp}(k,\ell)\rightarrow \mathrm{Mp}(k|\ell)$ be the
projection. On the upshot, we have:
\begin{lem}
There is a commutative diagram of morphisms of Lie groups
\begin{equation}
\begin{array}{ccccc}
    & & \kappa _1& & \\
   & \mathrm{Mp}(k:\ell) & \longrightarrow & \mathrm{ML}(k)\times \mathrm{Mp}(\ell -k) \\
   &   &  &  & \\
 \kappa_2&\downarrow &\stackrel{\tilde o}{\swarrow}& \downarrow & \mu\\
   &   & \theta &  & \\
   &\mathrm{Mp}(k,\ell)  & \longrightarrow & \mathrm{Mp}(k|\ell)& \\
   &                     &           &                    & \\
  s &\downarrow&  &\downarrow  & \nu\\
   &   &\varpi&  & \\
   &\mathrm{Sp}(k,\ell)&\longrightarrow&\mathrm{GL}(k)\times \mathrm{Sp}(\ell -k)          &\\
\label{eqn:fibreprodLie}
\end{array}\end{equation}
where the upper square is a fibre product diagram, $\mu$ and $\nu$
are 2:1 coverings, and $\nu \circ \mu$ is the product of the 2:1
coverings $\mathrm{ML}(k)\rightarrow \mathrm{GL}(k)$ and
$\mathrm{Mp}(\ell-k)\rightarrow \mathrm{Sp}(\ell-k)$.
\label{lem:fibreprodLie}\end{lem}

Recall that for every integer $l\ge 1$ the metalinear group
$\mathrm{ML}(l)$ has four connected component. If
$r:\mathrm{ML}(l)\rightarrow \mathrm{GL}(l)$ is the projection,
the composition $\det \circ \,r\,:\,\mathrm{ML}(l)\rightarrow
\mathbb{R}$ admits a square root, that we shall denote by
$\sqrt{\det}:\mathrm{ML}(l)\rightarrow \mathbb{C}$. We shall
denote the four connected components of $\mathrm{ML}(l)$ by
$\mathrm{ML}(l)_{1}$, $\mathrm{ML}(l)_{-1}$, $\mathrm{ML}(l)_{i}$,
$\mathrm{ML}(l)_{-i}$, meaning that $\sqrt{\det}>0$ on
$\mathrm{ML}(l)_{1}$, that $\sqrt{\det}$ is positive imaginary on
$\mathrm{ML}(l)_{i}$, and so forth. Clearly, $\mathrm{ML}(l)_{1}$
is the identity component, and $r^{-1}\left (I_{\ell}\right
)=\{e,g_l\}$ where $e$ is the identity and $g_l\in
\mathrm{ML}(l)_{-1}$. Given the standard inclusion
$\mathrm{ML}(l)\hookrightarrow \mathrm{Mp}(l)$, we shall view
$g_l$ as sitting in $\mathrm{Mp}(l)$.

\begin{lem} \label{lem:ker}
$\ker (\mu)=\ker (\tilde o)=\{e,(g_k,g_{\ell-k})\}$, where
$e\in \mathrm{ML}(k)\times \mathrm{Mp}(\ell-k)$ is the unit.
\end{lem}

\noindent \textit{Proof.} Since $\mu$ is a 2:1 cover, it suffices
to show that $(g_k,g_{\ell-k})\in \ker (\tilde o)$. Since $s\circ
\tilde o\big ((g_k,g_{\ell-k})\big )\in T=\ker (\varpi)$, and $T$
has no elements of order 2, $s\circ \tilde o\big
((g_k,g_{\ell-k})\big )$ is the identity. Given that $\tilde o\big
((g_k,g_{\ell-k})\big )$ lies in $\mathrm{ML}(\ell)$, it suffices
to notice that $\sqrt{\det}\circ \tilde o\big
((g_k,g_{\ell-k})\big )=(-1)\cdot (-1)=1$.

\begin{cor} \label{cor:inj}
By restriction of $\tilde o$, we have an injective morphism of Lie
groups $\mathrm{GL}_+(k)\times \mathrm{Mp}(\ell-k)\hookrightarrow
\mathrm{Mp}(k,\ell)\subseteq \mathrm{Mp}(\ell)$.
\end{cor}

\bigskip

Let us identify the unitary group $\mathrm{U}(\ell)$ with the
maximal compact subgroup of $\mathrm{Sp}(\ell)$ consisting of all
$2\ell\times 2\ell$ symplectic matrices of the form
$$
\left(%
\begin{array}{cc}
  X & -Y \\
  Y & X \\
\end{array}%
\right) \text{ such that } X+iY\in \mathrm{U}(\ell).$$ Notice that
$T\cap \mathrm{U}(\ell)=\{I_{2\ell}\}$ and $Z\cap
\mathrm{U}(\ell)\cong \mathrm{O}(k)$.

Similarly, let us set $\mathrm{U}(k,\ell)=\mathrm{Sp}(k,\ell)\cap
\mathrm{U}(\ell)$ and
$\mathrm{U}_+(k,\ell)=\mathrm{Sp}_+(k,\ell)\cap \mathrm{U}(\ell)$,
so that
$$
\begin{array}{lll}
\mathrm{U}(k,\ell)&=& \left \{ \left (
\begin{array}{cccc}
  A & 0 & 0 & 0 \\
  0 & R & 0 & S \\
  0 & 0 & A & 0 \\
  0 & -S & 0 & R
\end{array}
\right ) \, :\, A\in \mathrm{O}(k),\,
\left(%
\begin{array}{cc}
  R & -S \\
  S & R
\end{array}%
\right)\in \mathrm{U}(\ell-k) \right \}\\
& \cong& \mathrm{O}(k)\times \mathrm{U}(\ell-k),\end{array} $$ and
$$
\mathrm{U}_+(k,\ell)\, \cong \, \mathrm{SO}(k)\times
\mathrm{U}(\ell-k)$$ is the identity component of
$\mathrm{U}(k,\ell)$.

Let $\mathrm{MU}(\ell)=s^{-1}(\mathrm{U}(\ell))\subseteq
\mathrm{Mp}(\ell)$. This is a connected compact subgroup of
$\mathrm{Mp}(\ell)$. Similarly, let
$\mathrm{MU}(k,\ell)=s^{-1}(\mathrm{U}(k,\ell))\subseteq
\mathrm{Mp}(k,\ell)$. If we let
$\mathrm{MO}(k)=s^{-1}(\mathrm{O}(\ell))$, the vertical maps in
the commutative diagram (\ref{eqn:fibreprodLie}) both reduce to
the composition of double covers
$$\mathrm{MO}(k)\times
\mathrm{MU}(\ell-k)\longrightarrow
\mathrm{MU}(k,\ell)\longrightarrow \mathrm{U}(k,\ell),$$ and the
analogue of Corollary \ref{cor:inj} is

\begin{cor}\label{cor:injunitary}
By restriction of $\tilde o$, we have an injective morphism of Lie
groups $\mathrm{SO}(k)\times \mathrm{MU}(\ell-k)\hookrightarrow
\mathrm{MU}_+(k,\ell)\subseteq \mathrm{MU}(\ell)$.
\end{cor}

For any $l\ge 1$, let $\mathcal{S}(\mathbb{R}^l)$ be the space of
smooth complex valued rapidly decaying functions on
$\mathbb{R}^l$, endowed with the bilinear pairing
$$(f,\,g)\,=:\, \int _{\mathbb{R}^l}\, f(x)\, g(x)\,
dx\,\,\,\,\,\,\,\,\,(f,g\in \mathcal{S}(\mathbb{R}^l)),$$ and the
$L^2$-Hermitian product $(f,g)_h=:(f,\overline g)$. The
metaplectic group acts unitarily on $\left
(\mathcal{S}(\mathbb{R}^l),\, (\,,\,)_h\right )$ under the
Segal-Shale-Weyl representation,
$$\upsilon _{\mathrm{SSW}}\, :\,\mathrm{Mp}(l)\longrightarrow
\mathrm{U}\left (\mathcal{S}(\mathbb{R}^l)\right ),$$ which plays
a crucial role in the symbolic calculus of Fourier-Hermite
distributions.

Given the inclusion of Lie groups of Corollary \ref{cor:inj} (or
Corollary \ref{cor:injunitary}), we may restricts the
Segal-Shale-Weyl representation of $\mathrm{Mp}(\ell)$ to
$\mathrm{GL}_+(k)\times \mathrm{Mp}(\ell-k)$. On the other hand,
$\mathrm{GL}_+(k)\times \mathrm{MU}(\ell-k)$ acts on
$\mathcal{S}(\mathbb{R}^{\ell-k})=\mathbb{C}\otimes _{\mathbb{C}}
\mathcal{S}(\mathbb{R}^{\ell-k})$ by the tensor product
$\sqrt{\det}\, \otimes \, \upsilon _{\mathrm{SSW}}$ of the
charachter $\sqrt{\det}$ on $\mathrm{GL}_+(k)\subseteq
\mathrm{ML}(k)$ and the Segal-Shale-Weyl representation of
$\mathrm{Mp}(\ell-k)$. Using the description of the metaplectic
representation in \cite{gs-ga} and \cite{bg}, one can check that

\begin{lem} The restriction map
$\mathcal{S}(\mathbb{R}^\ell)\longrightarrow
\mathcal{S}(\mathbb{R}^{\ell-k})$, $$f(x_1,\ldots,x_\ell)\mapsto
f_{\mathrm{res}}(x_{k+1},\ldots,x_\ell)=:f(0,\ldots,0,x_{k+1},\ldots,x_\ell)\,\,\,\,\,\,\,\,\,\,\,(f\in
\mathcal{S}(\mathbb{R}^\ell))$$ is equivariant, that is, it is a
morphism of $\mathrm{GL}_+(k)\times \mathrm{MU}(\ell-k)$-modules.
\label{lem:equiv}
\end{lem}

In the following, we shall use the concepts of metaplectic
manifold and metaplectic vector bundle, a manifold being
metaplectic if its tangent bundle is. If $(E,\Omega _E)$ is a
symplectic vector bundle of rank $2\ell$ over a manifold $N$, we
shall denote by $\mathrm{Bp}(E)\rightarrow N$ the principal
$\mathrm{Sp}(\ell)$-bundle of all symplectic frames in $E$. If the
symplectic vector bundle $(E,\Omega _E)$ is metaplectic, we shall
denote by $\widetilde{Bp}(E)$ the corresponding principal
$\mathrm{Mp}(\ell)$-bundle.

In particular, if $E$ admits a Lagrangian subbundle $L\subseteq
E$, then $E$ is symplectically equivalent to the vector bundle
$L\oplus L^*$, with its standard symplectic structure. The
structure group of $E$ then reduces to $\mathrm{GL}(\ell)\subseteq
\mathrm{Sp}(\ell)$, and $E$ is metaplectic if and only if $L$ is
metalinear. In particular, a cotangent bundle $T^*M$ is a
metaplectic manifold if and only the base manifold $M$ is
metalinear, and furthermore any Lagrangian submanifold of a
metaplectic manifold is metalinear.

Let us now assume, more generally, that $S\subseteq E$ is a
rank-$k$ isotropic vector subbundle. The subbundle
$\mathrm{Bp}(S,E)\subseteq \mathrm{Bp}(E)$ consisting of all
symplectic basis whose first $k$ vectors lie in $S$ is a principal
$\mathrm{Sp}(k,\ell)$-bundle over $N$. Let $\mathrm{BL}(S)$ denote
the principal $\mathrm{GL}(k)$-bundle over $N$ consisting of all
linear frames in $S$. Let $S^{\perp_\Omega}\subseteq E$ be the
symplectic annihilator of $S$. We have $\mathrm{BL}(S)\times
\mathrm{Bp}(S^{\perp_\Omega}/S)=\mathrm{Bp}(S,E)/T$, the
projection being the bundle map
\begin{equation}\label{eqn:bundlemorph}(e_1,\ldots,e_\ell,f_1,\ldots
f_\ell)\mapsto \big
((e_1,\ldots,e_\ell),([e_{k+1}],\ldots,[e_\ell],[f_{k+1}],\ldots,[f_\ell])\big
),\end{equation} obviously equivariant with respect to the
morphism of Lie groups (\ref{eqn:liegrmorph}).

Now suppose, in addition, that $S$ is orientable. Let
$\mathrm{BL}_+(S)\subseteq \mathrm{BL}(S)$ be the principal
$\mathrm{GL}_+(k)$-bundle of all oriented frames in $S$. Let
$\mathrm{Bp}_+(S,E)\subseteq \mathrm{Bp}(S,E)$ be the subbundle
consisting of all symplectic basis whose first $k$ vectors form an
\textit{oriented} basis of $S$. Then $\mathrm{Bp}_+(S,E)$ is a
principal $\mathrm{Sp}(k,\ell)_+$-bundle over $N$, and the
projection $\mathrm{Bp}_+(S,E)\longrightarrow
\mathrm{BL}_+(S)\times \mathrm{Bp}(S^{\perp_\Omega}/S)$ is
equivariant with respect to the morphism of Lie groups $\varpi _+$
in (\ref{eqn:liegrmorph+}).

\begin{prop} Let $(E,\Omega _E)$ be a symplectic vector bundle,
and let $I\subseteq E$ be an oriented rank-$k$ isotropic vector
subbundle. Let $I^{\perp _\Omega}\subseteq E$ denote the
symplectic annihilator of $S$ in $E$. Then there is a natural
bijection between the set of equivalence classes of metaplectic
structures on the symplectic vector bundle
$N_I=:I^{\perp_\Omega}/I$ and the set of of equivalence classes of
metaplectic structures on $E$. \label{prop:metaplbiject}
\end{prop}

\noindent \textit{Proof.} In one direction, suppose given a
metaplectic structure on $E$, that we describe by the following
equivariant commutative diagram of principal bundles and double
covers:
$$
 \begin{array}{ccc}
    \mathrm{Mp}(\ell)\times \mathrm{Mp}(E) & \longrightarrow &\mathrm{Mp}(E) \\
     & & \\
\downarrow& & \downarrow\\
     & & \\
   \mathrm{Sp}(\ell)\times \mathrm{Bp}(E) & \longrightarrow
&\mathrm{Bp}(E).
  \end{array}
$$
By restriction, we obtain the other
$$
 \begin{array}{ccc}
    \mathrm{Mp}_+(k,\ell)\times \mathrm{Mp}_+(I,E) & \longrightarrow &\mathrm{Mp}_+(I,E) \\
     & & \\
\downarrow& & \downarrow\\
     & & \\
   \mathrm{Sp}_+(k,\ell)\times \mathrm{Bp}_+(I,E) & \longrightarrow
&\mathrm{Bp}_+(I,E).
  \end{array}
$$
Given the projection
$$\mathrm{Bp}_+(I,E)\, \longrightarrow \, \mathrm{Bp}_+(I,E)/Z=\mathrm{Bp}(N_I),$$
Lemma \ref{lem:metapl} now shows that
$\mathrm{Mp}(N_I)=:\mathrm{Mp}_+(I,E)/Z$ is a metaplectic
structure on $N_I$.

In the opposite direction, suppose given a metaplectic structure
on $N_I$. Let $\mathcal{J}(E,\Omega_E)$ be the contractible space
of all complex structures on $E$ compatible with $\Omega _E$
\cite{ms}, and let us fix $J\in \mathcal{J}(E,\Omega_E)$. Thus $E$
inherits a unitary structure, say $h=g+i\,\Omega$ where
$g(v,w)=\Omega (v,Jw)$ and, by restriction of $g$, $S$ is a
Euclidean vector bundle. Let $I^{\perp_g}$ and $I^{\perp _h}$ be
the Euclidean and Hermitian orthocomplements of $I$ in $E$. The
intersection $I^{\perp_h}=I^{\perp _g}\cap I^{\perp _\Omega}$ is a
complex vector subbundle of $E$, linearly symplectomorphic to the
vector bundle $N_I$. Having this identification in mind, let
$\mathrm{BU}(N_I)$ be the bundle of all unitary frames in $N_I$.
The metaplectic structure then yields an equivariant commutative
diagram
$$
\begin{array}{ccc}
    \mathrm{MU}(k)\times \widetilde{\mathrm{BU}}(N_I) &\longrightarrow & \widetilde{\mathrm{BU}}(N_I)\\
& & \\
\downarrow & & \downarrow\\
& & \\
 \mathrm{U}(k)\times
\mathrm{BU}(N_I) &\longrightarrow & \mathrm{BU}(N_I).
   \end{array}
$$
Let $\mathrm{BO}_+(I)\subseteq \mathrm{BL}_+(I)$ be the subbundle
of all oriented orthonormal frames of $I$. Let
$\mathrm{BU}_+(I,E)=\mathrm{BU}(E)\cap \mathrm{BL}_+(I,E)$. We
have an isomorphism of principal $\mathrm{U}_+(k,\ell)$-bundles
$$\mathrm{Bp}(E)\supseteq \mathrm{BU}_+(I,E)\,\cong
\,\mathrm{BO}_+(I)\times_N \mathrm{BU}(N_I).$$ Now
$\widetilde{\mathrm{BU}}_+(I,E)=\mathrm{BO}_+(I)\times_N
\widetilde{\mathrm{BU}}(N_I)$ is an equivariant double cover of
$\mathrm{BU}_+(I,E)$ with respect to the double cover of structure
groups $$\mathrm{SO}(k)\times \mathrm{MU}(\ell-k)\rightarrow
\mathrm{U}_+(k,\ell),$$ and we obtain a metaplectic structure of
$E$ by group extension in view of Corollary \ref{cor:injunitary}.

It is clear that these two correspondences are one the inverse of
the other.

\bigskip

\begin{cor} Let $E$ be a metaplectic vector bundle of
rank $2\ell$. Let $S\subseteq E$ be an oriented rank-$k$ isotropic
subbundle. Then the structure group of $E$ reduces to
$\mathrm{GL}_+(k)\times \mathrm{Mp}(\ell-k)\subseteq
\mathrm{Mp}(\ell)$.\label{lem:reducedstrgrp} \end{cor}

If $E$ is metaplectic vector bundle of rank $2\ell$, given the
Segal-Shale-Weyl representation we may form the associated
infinite dimensional vector bundle $$
\mathcal{S}(E)=:\widetilde{\mathrm{Bp}}(E)\times
_{\mathrm{Mp}(\ell)}\mathcal{S}(\mathbb{R}^{\ell}).$$ If
$I\subseteq E$ is an oriented rank-$k$ isotropic subbundle, in
view of the induced metaplectic structure on $N_I$ we may
similarly form infinite dimensional vector bundle $$
\mathcal{S}(N_I)=:\widetilde{\mathrm{Bp}}(N_I)\times
_{\mathrm{Mp}(\ell-k)}\mathcal{S}(\mathbb{R}^{\ell-k}).$$ We have:

\begin{cor} \label{cor:surjmorphi}
Let $E$ be a metaplectic vector bundle of rank $2\ell$, and let
$I\subseteq E$ be an oriented rank-$k$ isotropic subbundle. Then
the restriction map of Lemma \ref{lem:equiv} extends to a
surjective morphism of vector bundles
$$
\Phi _I:\mathcal{S}(E)\, \longrightarrow \,\bigwedge
^{-1/2}(I)\otimes _{\mathbb{C}}\mathcal{S}(N_I).
$$
\end{cor}

By the bilinear pairing $(\,,\,)$, we may view in the standard
manner $\mathcal{S}(\mathbb{R}^\ell)$ as a subspace of the space
$\mathcal{S}'(\mathbb{R}^\ell)$ of tempered distribution, which is
invariant under the dual representation of $\mathrm{Mp}(\ell)$ on
$\mathcal{S}'(\mathbb{R}^\ell)$. If $E$ is as above a rank-$2\ell$
metaplectic vector bundle, let $\mathcal{S}_{\mathrm{dual}}(E)$
denote the vector associated to the metaplectic structure and the
dual action of $\mathrm{Mp}(\ell)$ on
$\mathcal{S}(\mathbb{R}^\ell)$. Since the metaplectic
representation is unitary, we have an isomorphism
$$
\mathcal{S}_{\mathrm{dual}}(E)\,\cong\, \overline{\mathcal{S}(E)}.
$$
Passing to conjugate bundles and inserting this isomorphism, we
obtain from Corollary \ref{cor:surjmorphi} a surjective morphism
of vector bundles \cite{bg}
\begin{equation}\label{eqn:conjugatesurj}
\mathcal{S}_{\mathrm{dual}}(E)\, \longrightarrow
\,\overline{\bigwedge} ^{-1/2}(I)\otimes
_{\mathbb{C}}\mathcal{S}_{\mathrm{dual}}(N_I).
\end{equation}

Consider now the inclusion of Lie groups $$
\mathrm{SO}(\ell)\subseteq \mathrm{GL}_+(\ell)\subseteq
\mathrm{ML}(\ell)\subseteq \mathrm{Mp}(\ell).$$ The function
$e^{-\|X\|^2}=e^{-\sum _ix_i^2}\in \mathcal{S}(\mathbb{R}^\ell)$
is a fixed point for the restriction of the Segal-Shale-Weyl
representation to the subgroup $\mathrm{SO}(\ell)$. In the special
case where $E$ admits an oriented Lagrangian subbundle $L$, as we
have mentioned the structure group reduces to
$\mathrm{GL}_+(\ell)$. The choice of a compatible complex
structure $J\in \mathcal{J}(E,\Omega_E)$ further reduces it to
$\mathrm{SO}(\ell)$, the corresponding principal
$\mathrm{SO}(\ell)$-bundle $\mathrm{P}_{L,J}$ being the bundle of
all frames of the type
$$(e_1,\ldots,e_\ell,J_p(e_1),\ldots,J_p(e_\ell)),$$ where
$(e_1,\ldots,e_\ell)$ is an oriented orthonormal frame for $L(p)$,
$p\in N$. Thus:

\begin{lem} Suppose that the symplectic bundle $(E,\Omega _E)$ has
an oriented Lagrangian subbundle $L\subseteq E$. Then to any $J\in
\mathcal{J}(E,\Omega_E)$ there is associated a natural nowhere
vanishing section $\sigma_{L,J}$ of $\mathcal{S}(E)$, which may be
described as the constant $\mathcal{S}(\mathbb{R}^\ell)$-valued
function on $\mathrm{P}_{L,J}$ equal to $e^{-\|x\|^2}$. Its image
in $\Omega ^{-1/2}(L)$ corresponds to the constant function $1$ on
$\mathrm{P}_{L,J}$.\label{lem:lagsubbundlesect}\end{lem}

In a related vein, the theory of \cite{bg} also shows that the
function $e^{-\|X\|^2}$ is a joint eigenvector for the metaplectic
action of $ \mathrm{MU}(\ell)=:s^{-1}\big (U(\ell)\big )\subseteq
\mathrm{Mp}(\ell)$, i.e. there is a unitary charachter
$c:\mathrm{MU}(\ell)\rightarrow \mathrm{U}(1)$ such that
$$\upsilon _{\mathrm{SSW}}(A)\left (e^{-\|X\|^2}\right )=c(A)\cdot
e^{-\|X\|^2}\,\,\,\,\,\,\,\,\,\,(A\in \mathrm{MU}(\ell).)$$ Thus,
the function $e^{-\|X\|^2-\|Y\|^2}$ is a fixed point for the
tensor product action of $\mathrm{MU}(\ell)$ on
$\mathcal{S}(\mathbb{R}^\ell)\otimes
\mathcal{S}_{\mathrm{dual}}(\mathbb{R}^\ell)\cong
\mathcal{S}(\mathbb{R}^\ell)\otimes
\overline{\mathcal{S}(\mathbb{R}^\ell)}$. It represents the
orthogonal projection of $\mathcal{S}(\mathbb{R}^\ell)$ onto the
subspace $\mathrm{span}\left \{e^{-\|X\|^2}\right \}$.

Now the choice of a compatible complex structure $J\in
\mathcal{J}(E,\Omega_E)$ reduces the structure group of $E$ to
$U(\ell)$; let $\mathrm{Bu}(E,J)\subseteq \mathrm{Bl}(E)$ be the
principal $\mathrm{U}(\ell)$-bundle of all unitary frames in $E$
(for $(\Omega _E,J)$). If $(E,\Omega_E)$ in addition is
metaplectic, with metaplectic strucure
$\widetilde{\mathrm{Bp}}(E)\stackrel{\varpi}{\rightarrow}\mathrm{Bp}(E)$,
the inverse image $\widetilde{\mathrm{Bu}}(E,J)=:\varpi ^{-1}\left
(\mathrm{Bu}(E,J)\right)$ is a principal
$\mathrm{MU}(\ell)$-bundle. Thus, we may globalize the previous
construction to obtain:

\begin{lem} Let $(E,\Omega _E)$ be a metaplectic vector bundle.
Then to each $J\in \mathcal{J}(E,\Omega_E)$ there is associated a
line subbundle $L_J\subseteq \mathcal{S}(E)$, which in the
trivialization of $\mathcal{S}(E)$ offered by any element of
$\widetilde{\mathrm{Bu}}(E,J)$ is the line spanned by
$e^{-\|X\|^2}$. In the same trivialization, the orthogonal
projector $\mathcal{S}(E)\rightarrow L_J$, as a section of the
bundle of linear endomorphisms of $\mathcal{S}(E)$, is represented
by the function $e^{-\|X\|^2-\|Y\|^2}$.
\end{lem}

When operating on symbols of Szeg\"{o} kernels and related
distributions, a naturally occurring case is that of symplectic
vector bundles of the form $E^+\oplus E^-$, as in the following
Lemma (that we simply state):

\begin{lem}

\begin{description}
  \item[i)] Let $E^+=(E,\Omega _E)$ and $F^+=(F,\Omega_F)$ be metaplectic
vector bundles on a manifold $N$. Then there is a naturally
induced metaplectic structure on their direct sum $E^+\oplus
F^+=(E\oplus F,\Omega _E\oplus \Omega _F)$, and a natural
isomorphism
$$\mathcal{S}(E^+\oplus F^+)\cong \mathcal{S}(E)\otimes
\mathcal{S}(F).$$
  \item[ii)] Let $E^+=(E,\Omega _E)$ be a metaplectic vector
bundle. Then there is a naturally induced metaplectic structure on
its opposite $E^-=:(E,-\Omega _E)$, and a natural isomorphism
$$\mathcal{S}(E^-)\cong \mathcal{S}_{\mathrm{dual}}(E^+).$$
  \item[iii)] In particular, if $E^+=(E,\Omega _E)$ is a metaplectic vector
bundle then there is a naturally induced metaplectic structure on
$E^+\oplus E^-$, and a natural isomorphism
$$\mathcal{S}(E^+\oplus E^-)\cong \mathrm{End}_{\mathcal{HS}}(\mathcal{S}(E)),$$
where the latter denotes the vector bundle of Hilbert-Schmidt
linear endomorphisms of $\mathcal{S}(E)$.
\end{description}\label{lem:usefulfacts}
\end{lem}

If $S=E^+\oplus E^-$, for any $J\in \mathcal{J}(E,\Omega _E)$ we
have $$J_S=:J\oplus (-J)\in \mathcal{J}\left (S,\Omega_E\oplus
(-\Omega _E)\right ).$$ Any unitary frame $(e_1,\ldots,e_\ell)\in
\mathrm{Bu}(E,J)=\mathrm{Bu}(E,-J)$ extends to a unitary frame
$(e_1,\ldots,e_\ell,e_1\ldots,e_\ell)\in \mathrm{Bu}(N,J_N)$.
Passing to metaplectic double covers,
$\widetilde{\mathrm{Bu}(E,J)}\subseteq
\widetilde{\mathrm{Bu}}(S,J_S)$. Summing up, we have:

\begin{cor} \label{cor:sympdirectsumpm}
Suppose that the symplectic vector bundle $(E,\Omega_E)$ is
metaplectic. Set $S=:E^+\oplus E^-$. Then $\mathcal{S}(S)\cong
\mathrm{End}_{\mathcal{HS}}(\mathcal{S}(E))$ has a distinguished
section $\sigma _J$ for every compatible complex structure $J\in
\mathcal{J}(E,\Omega _E)$. In any appropriate trivialization (in
the sense above), this is represented by the function
$e^{-\|X\|^2-\|Y\|^2}\in
\mathcal{S}(\mathbb{R}^{2\ell})$.\end{cor}

We now look at the image of the distinguished section in Corollary
\ref{cor:sympdirectsumpm} under the morphism of vector bundles in
Corollary \ref{cor:surjmorphi}. The proof of the following is a
case by case application of the local form in Lemma
\ref{lem:equiv} of the vector bundle morphism of Corollary
\ref{cor:surjmorphi}, and is left to the reader:

\begin{cor}\label{cor:imageEpmsurjmorph}
In the hypothesis of Corollary \ref{cor:sympdirectsumpm}, suppose
that $L\subseteq E$ is an oriented isotropic subbundle, and
consider the isotropic subbundle $\mathrm{diag}(L)=\{(l,l):l\in
L\}\subseteq S$. Fix $J\in \mathcal{J}(E,\Omega_E)$, and let $h$
be the Hermitian structure on $E$ associated to the compatible
pair $(\Omega _E,J)$. Let $g=\Re (h)$ be the associated Riemannian
metric. Let $L^\perp \subseteq E$ be the symplectic annihilator of
$L$ in $(E,\Omega _E)$, and let $L^0$ be the Riemannian
orthocomplement of $L$ in $(E,g)$. Then:
\begin{description}
  \item[i):] The intersection $L^\perp \cap L^0$ is a complex vector subbundle of
$(E,\Omega_E)$, whence a symplectic subbundle of $(E,\Omega _E)$;
  \item[ii):] The symplectic normal bundle $N_L=L^\perp/L$ of $L$
in $E^+$ is (naturally) symplectically isomorphic to $L^\perp \cap
L^0$ (we shall henceforth not distinguish between $N_L$ and
$L^\perp\cap L^0$), and thus has a naturally induced compatible
complex structure $J_{N_L}$;
  \item[iii):] Let us endow $S=E^+\oplus E^-$ with the compatible complex structure $J\oplus (-J)$.
Then the symplectic normal bundle $N_{\mathrm{diag(L)}}$ of
$\mathrm{diag}(L)$ in $S$ is (naturally isomorphic to) the direct
sum of vector subbundles
\begin{equation}\label{eqn:sympnormdiagonal}
  \begin{array}{ccl}
    N_{\mathrm{diag(L)}}& \cong & \big (\{(l,-l):l\in L\}\oplus \{(Jl,Jl):l\in L\}\big )\\
     & &\oplus \left (N_L^+\oplus N_L^-\right )\\
&=&(L_r\oplus L_i)\oplus \left (N_L^+\oplus N_L^-\right )\\
&=&L_{\mathbb{C}}\oplus \left (N_L^+\oplus N_L^-\right ).
  \end{array}
\end{equation}
Here $L_r$ and $L_i$, defined by the second identity, are oriented
Lagrangian subbundles of the complex vector subbundle
$L_{\mathbb{C}}=:L_r\oplus L_i\subseteq S$. Since $L_r\cong L$ is
oriented, $L_{\mathbb{C}}$ is metaplectic.
\item[iv):] Let $\Phi
_{\mathrm{diag}(L)}:\mathcal{S}(S)\rightarrow \bigwedge
^{-1/2}(L)\otimes \mathcal{S}\left (N_{\mathrm{diag}(L)}\right )$
be the vector bundle morphism introduced in Corollary
\ref{cor:surjmorphi}. Let $\sigma _J$ be the distinguished section
of $\mathcal{S}(S)$ associated to the complex structure $J$ as in
Corollary \ref{cor:sympdirectsumpm}. Then
$$
\Phi _{\mathrm{diag}(L)}\left ( \sigma _J\right
)=\mathrm{vol}_L^{-1/2}\otimes \sigma _{L_r}\otimes \sigma
_{J_{N_L}}.
$$
Here $\mathrm{vol}_L^{-1/2}$ is the section of $\bigwedge
^{-1/2}(I)$ taking value one on oriented orthonormal basis of $I$,
$\sigma _{L_r}$ is the section of $\mathcal{S}(L_{\mathbb{C}})$
associated to the Lagragian subbundle $L_r$ according to Lemma
\ref{lem:lagsubbundlesect}, and $\sigma _{J_{N_L}}$ is the section
of $\mathcal{S}(N_L^+\oplus N_L^-)$ associated to the complex
structure $J_{N_L}$ on $N_L$, according to Corollary
\ref{cor:sympdirectsumpm}.
\item[v):] The symplectic normal bundle of $L_r$ in $N_{\mathrm{diag}(L)}$
is isomorphic as a unitary vector bundle to $N_L^+\oplus N_L^-$.
Let $\Phi _{L_r}:\mathcal{S}(N_{\mathrm{diag}(L)})\rightarrow
\bigwedge ^{-1/2}(L_r)\otimes \mathcal{S}\left (N_{L_r}\right )$
be the vector bundle morphism from Corollary \ref{cor:surjmorphi}.
Then $$ \Phi _{L_r}\left ( \sigma _{L_r}\otimes \sigma
_{J_{N_L}}\right )=\mathrm{vol}_{L_r}^{-1/2}\otimes  \sigma
_{J_{N_L}}.
$$
\end{description}
\end{cor}

\section{The geometry of the symbol calculus}

We need to recall some basic constructions from \cite{bg}, at
places rephrasing them in terms of the principal bundles involved
in our constructions. Let $A$ and $B$ be $\mathcal{C}^\infty$
manifolds. There is a natural symplectomorphism $T^*(A\times
B)\cong T^*(A)\times T^*(B)$ (cotangent bundles are implicitly
endowed with their canonical symplectic structures), which will be
implicit throughout. Suppose that $\Gamma '\subseteq T^*(A\times
B)\setminus \{0\}$ is a closed Lagrangian conic submanifold. Let
$\Gamma \subseteq T^*(A)\times T^*(B)\{0\}$ be the associated
canonical relation, defined as the image of $\Gamma'$ under the
involution $\big ((a,\eta),(b,\vartheta)\big)\mapsto
((a,\eta),(b,-\vartheta)\big)$. Suppose that $\Sigma \subseteq
T^*(B)\setminus \{0\}$ is a closed isotropic conic submanifold.
Let us form the fibre diagram
\begin{eqnarray*}
  \begin{array}{ccccc}
   & & \rho & &\\
   &F&\longrightarrow&\Gamma& \\
   & &               &    &   \\
\varrho& \downarrow& & \downarrow &\gamma \\
   &  &      \iota        &    &   \\
   &\Sigma&\longrightarrow&T^*(B)\setminus \{0\}.&
   \end{array}
\end{eqnarray*}
Here $\iota$ is the inclusion, and $\gamma$ the projection. Thus,
$$
F=\left \{(\sigma, (\tau,\sigma))\, :\, \sigma \in \Sigma,\,
(\tau,\sigma)\in \Gamma\right \}.$$ We have a diffeomorphism
$$
F\, \cong \, q_\Gamma ^{-1}(\Sigma)\,=q^{-1}(\Sigma)\cap
\Gamma,\,\,\,\,\,\,\,\,(\sigma,(\tau,\sigma))\leftrightarrow
(\tau,\sigma),$$ so that $F$ maps naturally into $\Gamma$,
diffeomorphically onto its image. We shall implicitly think of $F$
as embedded in $\Gamma$, as the subset of all pairs
$(\tau,\sigma)\in \Gamma$ with $\sigma \in \Sigma$.

Let us assume that all the clean intersection and properness
hypothesis in Chapter 7 of \textit{loc. cit.} are satisfied (this
will be always the case in our situation). Given the projections
\begin{eqnarray*}
\begin{array}{cccccccccccc}
 &          & \Gamma &         & & & & & &T^*(A\times B) & & \\
 &\stackrel{p_\Gamma}{\swarrow}  &        &\stackrel{q_\Gamma}{\searrow} & & & & &\stackrel{p}{\swarrow} &       &\stackrel{q}{\searrow} &  \\
 T^*(A)&  &    & &T^*(B) & & &T^*(A) & &       & & T^*(B)
 \end{array}
 \end{eqnarray*}
we have a conic isotropic submanifold
\begin{eqnarray*}
\Gamma \circ \Sigma &=&\left \{\tau \in T^*(A)\setminus
\{0\}\,:\,\exists \,\sigma \in \Sigma\text{ such that }
(\tau,\sigma)\in \Gamma \right\}\\
&=& p\, \left (q^{-1}(\Sigma)\cap \Gamma\right )\\
& =&p_\Gamma (q_\Gamma ^{-1}(\Sigma)\big )\subseteq
T^*(A)\setminus \{0\}.
\end{eqnarray*}
Furthermore, the projection $p_F:F\rightarrow \Gamma\circ \Sigma$,
$(\sigma,(\tau,\sigma))\mapsto \tau$, is a fibration with compact
fibres.

Let us define a vector bundle $U_0$ on $F$ by
\begin{eqnarray*}
U_0(\tau,\sigma)&=&\left \{w\in T_\sigma \Sigma \,:\,(0,w)\in
T_{(\tau,\sigma)}\Gamma\right \}.
\end{eqnarray*}
Since $p_F$ is a submersion, if $w\in U_0(\tau,\sigma)$ then
$(0,w)\in T_{(\tau,\sigma)}F$ is tangent to the fibre of $p_F$
over $\tau$. In other words, $TF\supseteq U_0=\ker (dp_F)$ is the
vertical tangent bundle of $p_F$.

The functorial behaviour of the symbol of Fourier-Hermite
generalized half-forms is governed by a \textit{symbol map}, which
transforms symplectic spinors on $\Sigma$ into symplectic spinors
on $\Gamma \circ \Sigma$, and whose existence is the content of
Proposition 6.5 of \cite{bg}. More precisely, in the present
context this is a surjective morphism of vector bundles on $F$,
\begin{eqnarray*}
   \Psi_{\Sigma,\Gamma}: \rho ^*\left (\bigwedge ^{1/2}(\Gamma)\right )
\otimes _{\mathbb{C}}\varrho ^*\left ( \mathrm{Spin}(\Sigma)\right
) & \longrightarrow & \det (U_0^*)\otimes p_F^*\left (
\mathrm{Spin}(\Gamma \circ \Sigma)\right ).
    \end{eqnarray*}
In the sequel, we shall need to
reformulate its construction, in order to have an explicit
description of the symbol map in terms of the principal bundles
involved.

Since in our applications $A$, $B$, $\Gamma$, $\Sigma$, $F$, the
vector bundle $U_0$ above and the vector bundle $U$ introduced
below will all be orientable, we shall make this simplifying
assumption throughout. It will also simplify our exposition to
assume, as will be the case in our applications, that
$p_F:F\rightarrow \Gamma \circ \Sigma$ is a Riemannian submersion,
having therefore a natural connection, and that $TF$ has a natural
(oriented) complement in $\left .T\Gamma \right |_F$, denoted
$N_{F|\Gamma}$.

Given that $B$ is orientable, hence metalinear, $T^*B$ is
metaplectic. Since $\Sigma \subseteq T^*B\setminus \{0\}$ is an
oriented isotropic submanifold, its symplectic normal bundle in
$T^*B$, $N_\Sigma =(T\Sigma )^\perp /T\Sigma$, has an induced
metaplectic structure. The spinor bundle of $\Sigma$ is
$$
\mathrm{Spin}(\Sigma)=\bigwedge ^{1/2}(\Sigma)\otimes \mathcal{S}
(N_\Sigma).$$ The same considerations apply to the conic isotropic
submanifold $\Gamma \circ \Sigma \subseteq T^*A\setminus \{0\}$,
with associated spinor bundle
$$
\mathrm{Spin}(\Gamma \circ \Sigma)=\bigwedge ^{1/2}(\Gamma \circ
\Sigma)\otimes \mathcal{S} (N_{\Gamma \circ \Sigma}).$$

Let us define a second vector bundle $U_1\supseteq U_0$ on $F$ by
\begin{eqnarray*} U_1(\tau,\sigma)&=&\left \{w\in (T_\sigma
\Sigma)^\perp \,:\,(0,w)\in T_{(\tau,\sigma)}\Gamma\right \},
\end{eqnarray*}
for $(\tau,\sigma)\in F$.

Then $U_0$ and $U_1$ are isotropic vector subbundles of the
symplectic vector bundle $\varrho ^*\left (T(T^*B)\right )$, and
their quotient $U=U_1/U_0\subseteq \varrho ^*\left (N_\Sigma\right
)$ is an isotropic vector subbundle, which as mentioned we shall
assume oriented. By Corollary \ref{cor:surjmorphi}, therefore, the
structure group of the metaplectic bundle of $N_\Sigma$,
$\widetilde{\mathrm{Bp}}(N_\Sigma)$, reduces to
$$\mathrm{Bl}_+(U)\times \widetilde{\mathrm{Bp}}(U^\perp/U),$$
and there is furthermore a surjective morphism of vector bundles
associated to this reduced principal bundle,
$$
G_{\Sigma,U}:\mathcal{S}(N_\Sigma)\longrightarrow \bigwedge
^{-1/2}(U)\otimes \mathcal{S}\left (U^\perp /U\right ).$$ By
construction, $G_{\Sigma,U}$ is in the appropriate trivialization
the restriction map of Lemma \ref{lem:equiv}.

By Proposition 6.4 of \cite{bg}, $N_{\Gamma \circ \Sigma}\cong
U^\perp/U$ naturally, so that there is at any rate a morphism of
vector bundles on $F$ from $\rho ^*\left (\bigwedge
^{1/2}(\Gamma)\right ) \otimes _{\mathbb{C}}\varrho ^*\left (
\mathrm{Spin}(\Sigma)\right ) $ to
\begin{equation}
 \rho ^*\left (\bigwedge ^{1/2}(\Gamma)\right
)\otimes \varrho ^*\left (\bigwedge ^{1/2}(\Sigma)\right )\otimes
\bigwedge ^{-1/2}(U)\otimes p_F^*\left ( \mathcal{S}(N_{\Gamma
\circ \Sigma})\right ).
    \end{equation}
We shall therefore be done by exhibiting an isomorphism of line
bundles
\begin{equation}\label{eqn:isomolnbdl}
 \rho ^*\left (\bigwedge ^{1/2}(\Gamma)\right
)\otimes \varrho ^*\left (\bigwedge ^{1/2}(\Sigma)\right )\otimes
\bigwedge ^{-1/2}(U)\rightarrow \det (U_0^*)\otimes p_F^*\left
(\bigwedge ^{1/2}(\Gamma \circ \Sigma)\right ).
    \end{equation}
For the rest of this proof, we shall fix a point
$f=(\tau,\sigma)\in F$, and write for notational simplicity
$T\Sigma$, $T\Gamma$, and so forth for $T_fF$, $T_f\Gamma$, and so
forth. Similarly, $U_0$, $U_1$, $U$ will be the fibres at $f$ of
the corresponding vector bundles on $F$, and the same provision
will apply to the various principal bundles involved.

Given the assumptions discussed above, we have direct sum
decompositions
$$T\Gamma \cong TF\oplus N\cong  U_0\oplus T(\Gamma \circ \Sigma)\oplus
N,$$ whence we may reduce the principal bundle of the line bundles
on both sides of (\ref{eqn:isomolnbdl}) to the product
$$\mathrm{Bl}_+(U_0)\times \mathrm{Bl}_+(T(\Gamma \circ \Sigma))
\times \mathrm{Bl}_+(N)\times \mathrm{Bl}_+(\Sigma)\times
\mathrm{Bl}_+(U),$$ whose general element we shall denote by
$(e_{U_0},e_{\Gamma \circ \Sigma}, e_N,e_{\Sigma},e_U)$.

Let $U_1^\perp$ be the symplectic annihilator of $U_1$ in
$T(T^*B)$. Let $\gamma :T\Gamma \oplus T\Sigma \rightarrow U^\perp
_1$ be given by $\gamma ((v,w),w')=w-w'$. We then have a short
exact sequence (\cite{bg}, page 45):
$$0\,\longrightarrow\,TF\, \longrightarrow \,T\Gamma \oplus
T\Sigma \,\stackrel{\gamma}{\longrightarrow} \,U_1^\perp\,
\longrightarrow 0.$$ Under $\gamma$, the subspace $N\oplus T\Sigma
\subseteq T\Gamma\oplus T\Sigma$ maps isomorphically onto
$U_1^\perp \cong (T(T^*B)/U_1)^*$. Given a pair $(e_N,e_\Sigma)\in
\mathrm{Bl}_+(N)\times \mathrm{Bl}_+(\Sigma)$, we shall denote by
$f_{(e_N,e_\Sigma)}$ the oriented basis of $T(T^*B)/U_1$ dual to
the oriented basis $\gamma \big ((e_N,e_\Sigma)\big )$ of
$(T(T^*B)/U_1)^*$.

Let $\Omega _{T^*B,\mathrm{can}}$ denote the canonical symplectic
structure of $T^*B$, $\mathrm{vol}_{T^*B}=\Omega
_{T^*B,\mathrm{can}}^{\wedge b}$ the associated volume form, and
$\mathrm{vol}_{T^*B}^{1/2}=\sqrt{\Omega
_{T^*B,\mathrm{can}}^{\wedge b}}$ the associated half-form. If now
we are given $\vartheta _\Gamma \in \bigwedge ^{1/2}(\Gamma)$,
$\vartheta _\Sigma \in \bigwedge ^{1/2}(\Sigma)$, $\vartheta
_U^{-1}\in \bigwedge ^{-1/2}(U)$, the expression
\begin{equation} \label{eqn:halfformspiecedtogether}
\vartheta _\Gamma (e_{U_0},e_{\Gamma \circ \Sigma}, e_N)\cdot
\vartheta _\Sigma(e_{\Sigma})\cdot \vartheta ^{-1}_U(e_U)\cdot
\mathrm{vol}_{T^*B}^{1/2}(e_{U_0},e_U,f_{(e_N,e_\Sigma)})\end{equation}
only depends on $(e_{U_0},e_{\Gamma \circ \Sigma})\in
\mathrm{Bl}_+(U_0)\times \mathrm{Bl}_+(\Gamma \circ \Sigma)$.

We thus have:

\begin{prop}\label{prop:symbolmap}
Assume that $U_0$ is oriented. There the previous construction
defines a surjective morphism of vector bundles on $F$
  \begin{eqnarray*}
   \Psi_{\Sigma,\Gamma}: \rho ^*\left (\bigwedge ^{1/2}(\Gamma)\right )
\otimes _{\mathbb{C}}\varrho ^*\left ( \mathrm{Spin}(\Sigma)\right
) & \longrightarrow & \det (U_0^*)\otimes p_F^*\left (
\mathrm{Spin}(\Gamma \circ \Sigma)\right ).
    \end{eqnarray*}
\end{prop}

Our additional hypothesis that $U_0$ be orientable accounts for
the appearance of the determinant line bundle in place of the line
bundle of densities.

\section{Restricting and pushing forward $\tilde \Pi _X$.}

Let $L^*$ be the dual line bundle to $L$, endowed with the induced
Hermitian metric. Let $X\subseteq L^*$ be the unit circle bundle,
$H(X)\subseteq \mathcal{C}^\infty(X)$ be the Hardy space of
boundary values of holomorphic functions. The $S^1$-action induces
a decomposition $H(X)=\bigoplus _{k\in \mathbb{N}}H_k(X)$, and the
the $k$-th isotype $H_k(X)$ is canonically isomorphic to
$H^0(M,L^{\otimes k})$. Given that $\Omega$ is symplectic, the
connection form $\alpha$ is a contact structure on $X$. Given that
Hermitian metric $h$ on $L$ is $G$-invariant, the action of $G$ on
$L^*$ leaves $X$ invariant.

The K\"{a}hler manifold $(M,\Omega)$ and the contact manifold
$(X,\alpha)$ have natural volume forms $\mathrm{vol}_M=\Omega
^{\wedge n}$ and $\mathrm{vol}_X=\pi _X^*( \mathrm{vol}_M)\wedge
\alpha$, respectively. Let
$\mathrm{vol}^{(1/2)}_M=\sqrt{\mathrm{vol}_M}$ and
$\mathrm{vol}^{(1/2)}_X=\sqrt{\mathrm{vol}_X}$ denote the
respective associated half-forms. The volume form $\mathrm{vol}_X$
makes the space of smooth functions on $X$ into a prehilbert
vector space, unitarily isomorphic to the space $\Omega ^{1/2}(X)$
of smooth half-forms under the map $f\mapsto
f\,\mathrm{vol}^{(1/2)}_X$. This isomorphism will be implicit
throughout, and we shall accordingly view the Szeg\"{o} kernel
$\tilde \Pi _X$, rather than as a generalized density, as a
generalized half-form on $X\times X$: $\tilde \Pi _X\in
\mathcal{D}'_{1/2}(X\times X)$. Similarly, taking products we
obtain volume forms $\mathrm{vol}_{M\times M}$ and
$\mathrm{vol}_{X\times X}$ on $M\times M$ and $X\times X$,
respectively, with associated half-forms
$\mathrm{vol}^{(1/2)}_{M\times M}$ and
$\mathrm{vol}^{(1/2)}_{X\times X}$.

In general, if $f_1$ and $f_2$ are functions, or sections, or
half-forms, and so forth, on manifolds $A_1$ and $A_2$, we shall
denote by $f_1\boxtimes f_2=\pi _1^*(f_1)\otimes \pi _2^*(f_2)$
the corresponding object on the product $A_1\times A_2$ obtained
by pull-back under the projections $\pi _i:A_1\times
A_2\rightarrow A_i$ and tensor product. Thus, we may write the
Szeg\"{o} kernel as
$$\tilde \Pi _X \, =\, \sum _{k,j}\,s_j^{(k)}\boxtimes
\overline{s_j^{(k)}}\cdot \mathrm{vol}^{(1/2)}_{X}\boxtimes
\mathrm{vol}^{(1/2)}_{X},$$ where for every $k=1,2,\ldots$ $\left
\{s_j^{(k)}\right \}$ ($1\le j\le h^0(M,L^{\otimes k})$) is an
orthonormal basis for the $k$-th Fourier component $H(X)_k\cong
H^0(M,L^{\otimes k})$ of the Hardy space $H(X)$.

Let $M'=:\Phi ^{-1}(0)\subseteq M$. If $\pi _{X|M}:X\rightarrow M$
is the projection, let $X'=:\pi _{X|M}^{-1}(M')$ and denote by
$\jmath _{X'}:X' \hookrightarrow X$ be the inclusion. Under the
above assumptions, $p _{X'}:X'\rightarrow X_0=:X'/G$ is a
principal $G$-bundle. Clearly, $X_0$ is the unit circle bundle in
$L_0^*$ with the induced metric $h_0$. We shall now endow the
inclusion $\jmath _{X'}$ and the projection $p _{X'}$ with the
structure of morphisms of metalinear manifolds; by taking
products, this will also make the inclusion $\jmath _{X'\times
X'}:X'\times X'\rightarrow X\times X$ and the projection $\pi
_{X'}:X'\times X'\rightarrow X_0\times X_0$ into morphisms of
metalinear manifolds. We shall then apply the associated pull-back
and push-forward operations to the Szeg\"{o} kernel $\tilde \Pi
_X\in \mathcal{D}'_{1/2}(X\times X)$. By using the microlocal
theory of \cite{bg}, we shall relate the result to the Szeg\"{o}
kernel of the symplectic quotient $X_0$, $\Pi _{X_0}\in
\mathcal{D}'_{1/2}(X_0\times X_0)$.

To this end, let us fix an orientation on $G$, and thus on its Lie
algebra, $\frak{g}$. In particular, this makes $M'$ and $X'$ into
oriented Riemannian manifolds. As discussed in section 2, there is
then a natural way to restrict a half-form on $M$ (or $X$) to a
half-form on $M'$ (respectively, on $X'$), namely by setting
$\iota ^*(\eta _M)=\eta _{M'}$ (respectively, $\jmath ^*(\eta
_{X})=\eta _{X'}$). This makes $\iota$ and $\jmath$ into morphisms
of metalinear manifolds, $\tilde \iota$ and $\tilde \jmath$.

Taking products, we clearly obtain restriction maps $\iota
^*\otimes \iota ^*$ and $\jmath ^*\otimes \jmath ^*$ on half-form
bundles and $(\iota \times \iota )^*:\Omega ^{1/2}(M\times
M)\rightarrow \Omega ^{1/2}(M'\times M')$, $(\jmath \times \jmath
)^*:\Omega ^{1/2}(X\times X)\rightarrow \Omega ^{1/2}(X'\times
X')$ on smooth global half-forms.

Let $\Gamma _{\jmath\times \jmath}\subseteq T^*(X\times X)\times
T^*(X'\times X')$ be the canonical relation associated to
$\jmath\times \jmath$. The morphism of metalinear manifolds
$\tilde \jmath \times \tilde \jmath$ is equivalent to the
assignment of the half-form $\mathrm{vol}_{\Gamma _{\jmath\times
\jmath}}^{1/2}$ described in (\ref{eqn:halfformmeta}), with
$f=\jmath\times \jmath$, $P=X'\times X'$ and $Q=X\times X$.

By the same token, the principal $G$-bundle $p_{X'}:X'\rightarrow
X_0=X'/G$ is also a morphism of metalinear manifolds in the
natural manner: $p_{X'}^*(\eta _{X_0})=\eta _{X'}$. Let
$p_{X'*}:\Omega ^{1/2}(X')\rightarrow \Omega ^{1/2}(X_0)$ denote
the resulting push-forward operator. For any smooth function $f$
on $X'$, let us denote by $f^G$ the $G$-invariant component of
$f$, when decomposed over the irreducible representations of $G$.
Thus $f^G$ may be viewed implicitly as a function on $X_0=X'/G$.
Furthermore, let us write $V_{\mathrm{eff}}$ for the
\textit{effective potential} of the action, that is,
$V_{\mathrm{eff}}(\overline x)$ is the volume of the fibre $p_{X'}
^{-1}(\overline x)$, $\overline x\in X_0$ \cite{burnsg}. By Lemma
\ref{lem:propersubmersion}, if $\overline x\in X_0$ then
\begin{eqnarray*}\begin{array}{lcr}
p_{X'*}\left (f\,\mathrm{vol}^{(1/2)}_{V(X'/X_0)}\otimes
\mathrm{vol}^{(1/2)}_{H(X'/X_0)}\right )(\overline x)& & \\
&&\\
=\left ( \int _{p_{X'} ^{-1}(\overline x)}\,f \cdot
\mathrm{vol}^{(1/2)}_{V(X'/X_0)}\cdot
\overline{\mathrm{vol}^{(1/2)}_{V(X'/X_0)}}\right )\cdot
\mathrm{vol}_{X_0}^{(1/2)}&&\\
&=&f^G(\overline x)\cdot V_{\mathrm{eff}}(\overline x)\cdot
\mathrm{vol}_{X_0}^{(1/2)} \end{array}\end{eqnarray*}


Taking products, the principal $G\times G$-bundle $p_{X'\times
X'}=p_{X'}\times p_{X'}:X'\times X'\rightarrow X_0\times X_0$ can
be similarly made into a morphism of metalinear manifolds, with an
attached push-forward operation $p_{X'\times X'}:\Omega
^{1/2}(X'\times X')\rightarrow \Omega ^{1/2}(X\times X)$. If
$f,\,g$ are smooth functions on $X'$ and $(f\boxtimes
g)(x,y)=:f(x)\cdot g(y)$ ($x,y\in X'$), then $(f\boxtimes
g)^{G\times G}=f^G\boxtimes g^G$. Furthermore, the effective
potential for the product is
$V'_{\mathrm{eff}}=V_{\mathrm{eff}}\boxtimes V_{\mathrm{eff}}$,
that is, $V_{\mathrm{eff}}'(\overline x,\overline
y)=V_{\mathrm{eff}}(\overline x)\cdot V_{\mathrm{eff}}(\overline
y)$ ($\overline x,\, \overline y\in X_0$).

Thus,
\begin{lem}\label{lem:pushforward}
$\mathrm{vol}^{(1/2)}_{X'\times
X'}=\mathrm{vol}^{(1/2)}_{X'}\boxtimes \mathrm{vol}^{(1/2)}_{X'}$,
and
\begin{eqnarray*}
p_{X'\times X'*}\left (f\boxtimes
g\,\mathrm{vol}^{(1/2)}_{X'\times X'}\right ) &=& f^G\boxtimes
g^G\cdot V_{\mathrm{eff}}\boxtimes V_{\mathrm{eff}}\cdot
\mathrm{vol}_{X_0}^{(1/2)}\boxtimes
\mathrm{vol}_{X_0}^{(1/2)},\end{eqnarray*} for any pair of smooth
functions $f$ and $g$ on $X_0$. \end{lem}

\noindent \textit{Proof.} This follows from the equalities
\begin{eqnarray*}
p_{X'\times X'*}\left (f\boxtimes
g\,\mathrm{vol}^{(1/2)}_{X'\times X'}\right ) &=& p_{X'\times
X'*}\left ((f\cdot \mathrm{vol}^{(1/2)}_{X'})\boxtimes
(g\,\mathrm{vol}^{(1/2)}_{X'})\right )\\
&=& (f^G\,V_{\mathrm{eff}}\mathrm{vol}_{X_0 }^{(1/2)})\boxtimes
(g^G\, V_{\mathrm{eff}}\cdot \mathrm{vol}_{X_0}^{(1/2)}).
\end{eqnarray*}

\bigskip

Since $p_{X'\times X'}$ is proper, the push-forward $p_{X'\times
X'*}$ extends to a continuous linear map of Fr\'{e}chet vector
spaces
$$p_{X'\times X'*}:\mathcal{D}'_{1/2}(X'\times X')\,
\longrightarrow\,\mathcal{D}'_{1/2}(X_0\times X_0).$$
Let furthermore $\Sigma \subseteq T^*(X\times X)\setminus \{0\}$ be
the wave front set of the Szeg\"{o} kernel, that is,
$$\Sigma
=\left \{(x,x,r\alpha _x,-r\alpha _x)\,:\,x\in X,r>0\right \}.$$
Let $\mathcal{D}_{1/2}'(X\times X)_{\Sigma}\subseteq
\mathcal{D}_{1/2}'(X\times X)$ be the subspace of all generalized
half-forms on $X\times X$ having wave front contained in $\Sigma$, with the
appropriate topology. Let
\begin{equation*}
N_{\jmath\times \jmath} =\left \{(x_1,x_2,\eta _1,\eta_2): x_i\in
X',\, \eta _i=0\,\text{ on }T_{x_i}X'\subseteq T_{x_i}X \right \}
\end{equation*}
be the \textit{conormal bundle} of the embedding $\jmath\times
\jmath$. Since $\Sigma\cap N_{\jmath\times \jmath} =\emptyset$,
the pull-back $(\jmath \times \jmath)^*$ extends to a continuous
linear map of Fr\'{e}chet vector spaces (\cite{duist}, Proposition 1.3.3)
$$
(\jmath \times \jmath)^*:\,\mathcal{D}_{1/2}'(X\times X)_{\Sigma}\,
\longrightarrow \,\mathcal{D}_{1/2}'(X'\times X').$$
By composition, we thus have a well-defined continuous linear map
\begin{equation*}
\Upsilon =:p_{X'\times X'*}\circ (\jmath \times \jmath)^*:\,
\mathcal{D}_{1/2}'(X\times X)_{\Sigma}\, \longrightarrow
\,\mathcal{D}_{1/2}'(X_0\times X_0) .
\end{equation*}
Let us define, for $k\in \mathbb{Z}$,
\begin{eqnarray*}
  \mathcal{D}'_{1/2}(X_0\times X_0)_k &=& \left \{\varphi \in  \mathcal{D}'(X_0\times X_0):\right .\\
 &  & \left .u(e^{i\theta}x,y)=e^{ik\theta}u(x,y)\, \forall\, e^{i\theta}\in S^1,
 \forall x,y\in X\right \},
\end{eqnarray*}
and similarly on $\mathcal{D}'_{1/2}(X\times X)$ and
$\mathcal{D}'_{1/2}(X'\times X')$. Let
$T_k:\mathcal{D}_{1/2}'(X_0\times X_0)\,\longrightarrow
\,\mathcal{D}'_{1/2}(X_0\times X_0)_k$ be given by
\begin{eqnarray*}
  T_k(u) &=& \frac{1}{2\pi}\int _{0}^{2\pi}\, e^{-ik\theta}\,u(e^{i\theta}x,y)\,d\theta. \\
  \end{eqnarray*}
    Then
  \begin{equation*}
\mathcal{D}'_{1/2}(X_0\times X_0)\, \cong \, \bigoplus
_{k=-\infty}^{+\infty}\, \mathcal{D}'_{1/2}(X_0\times X_0)_k,
\end{equation*}
and $T_k$ represents the projection onto the $k$-th factor.

We now focus on the generalized half-form
$$
\tilde \Pi _{X|X_0}=:\Upsilon \left (\tilde \Pi _X\right )\in
\mathcal{D}'_{1/2}(X_0\times X_0).$$ Since $\jmath$ and $p_{X'}$
are $S^1$-equivariant maps, we have for every $k\in \mathbb{Z}$:
$$T_k\left (\Upsilon\left (\tilde \Pi _X\right )\right)\, =\,\Upsilon
\left (T_k(\tilde \Pi _X)\right ).$$ Clearly, $T_k(\tilde \Pi _X)$
is the distributional kernel of the orthogonal projector onto the
subspace $H(X)_k\cong H^0(M,L^{\otimes k})$.

For every $k=0,1,2,\cdots$ let $\left \{s_j^{(G,k)}\right \}_{1\le
j\le d_k}$ denote an orthonormal basis of the space of
$G$-invariant Hardy functions of level $k$, $H(X)_k^G\cong
H^0(M,L^{\otimes k})^G$. By Lemma \ref{lem:pushforward}, we have
\begin{equation}\label{eqn:localformpf}
\tilde \Pi
_{X|X_0}(x_0,y_0)=V_{\mathrm{eff}}(x_0)V_{\mathrm{eff}}(y_0)\,\sum
_{k=0}^{+\infty}\sum _{j=1}^{d_k}s_j^{(G,k)}(x')\otimes
\overline{s_j^{(G,k)}(y')},\end{equation} if $x_0,\,y_0\in X_0$
and $x',y'\in X'$ lie over $x_0,y_0$ respectively.

We shall now use the theory of \cite{bg} to describe $\tilde \Pi
_{X|X_0}$ as a Fourier-Hermite generalized half-form. Namely, let
us first recall that, as a Fourier-Hermite generalized half-form,
the Szeg\"{o} kernel satisfies $\tilde \Pi _X\in J^{1/2}(X\times
X,\Sigma)$. Let furthermore $\alpha _0$ be the connection 1-form
on $X_0$, and set
\begin{equation*}
 \Sigma _0=\left \{(x,x,r\alpha _{0\,x},-r\alpha _{0\,x}):x\in X_0,\,
r>0\right \}\subseteq T^*(X_0\times X_0)\setminus \{0\}.
\end{equation*}

\begin{lem}\label{lem:fourierhermite}
$\tilde \Pi _{X|X_0}\in J^{(1+g)/2}(X_0\times X_0,\Sigma _0)$.
\end{lem}

\noindent \textit{Proof.} Referring to chapter 7 of \cite{bg}, let
us consider the fibre diagram associated to $\Sigma \subseteq
T^*(X\times X)\setminus \{0\}$ and the conormal bundle $\Gamma
_{\jmath \times \jmath}'\subseteq T^*(X'\times X'\times X\times X)
\setminus \{0\}$ of the graph of $\jmath \times \jmath$:
\begin{equation}\label{eqn:fibdiag}
  \begin{array}{ccc}
    F_{\jmath \times \jmath} & \longrightarrow & \Gamma _{\jmath \times \jmath} \\
     &  & \\
    \downarrow & & \downarrow \\
     &   & \\
    \Sigma & \longrightarrow & T^*(X\times X),
  \end{array}
\end{equation}
where $\Gamma _{\jmath \times \jmath} $ denotes the image of
$\Gamma _{\jmath \times \jmath}' $ under sign reversal in the
first component of $T^*\big ((X'\times X')\times (X\times X)\big
)$. We have
\begin{equation}\label{defn:canreljperj}
\begin{array}{ccl}
\Gamma _{\jmath \times \jmath}&=&\left \{\left (x_1,x_2,\jmath
(x_1),\jmath (x_2),
(d_{x_1}\jmath)^t(\eta _1),(d_{x_2}\jmath)^t(\eta _2),\eta_1, \eta_2\right ):\right .\\
&&\left .x_i\in X',\, \eta_i\in T^*_{\jmath(x_i)}(X),\,
i=1,2.\right \}.
\end{array}\end{equation}
Let us set $\alpha '=:\jmath ^*(\alpha)$. Then
\begin{eqnarray*}
F_{\jmath \times \jmath}=\left \{\left (x,x,\jmath (x),\jmath (x),
r\alpha '_x,-r\alpha '_x,r\alpha _{\jmath (x)}, -r\alpha _{\jmath
(x)}\right ):\right . \left .x\in X',\, r>0.\right \}.
\end{eqnarray*}
Thus $F_{\jmath \times \jmath}$ is diffeomorphic to its projection
\begin{equation*}
 \Sigma '=\Gamma _{\jmath \times \jmath}\circ \Sigma=\left \{\left (x,x,
r\alpha '_x,-r\alpha '_x\right ):\right . \left .x\in X',\,
r>0.\right \}\subseteq T^*(X'\times X')\setminus \{0\}.
\end{equation*}
Therefore, the excess of the diagram (\ref{eqn:fibdiag}) is
\begin{eqnarray*}
  e _{\jmath \times \jmath}&=&\dim (F_{\jmath \times \jmath})+\dim T^*(X\times X)-\dim \left (\Gamma _{\jmath \times \jmath} '\right )
  -\dim (\Sigma) \\
  &=& g.
\end{eqnarray*}
Hence, given that all the clean intersection hypothesis of Theorem
9.1 of \cite{bg} are satisfied, we have $(\jmath \times \jmath
)^*(\tilde \Pi _X)\in J^{(1+g)/2}(X'\times X',\Sigma ')$.

Let us next consider the fibre diagram associated to $p_{X'}\times
p_{X'}$. The canonical relation is now
\begin{equation}\label{defn:canrelpperp}
\begin{array}{ccl}
\Gamma _{\jmath \times \jmath}&=&\left \{\left (x_1,x_2,p_{X'}
(x_1),p_{X'} (x_2), (d_{x_1}p_{X'})^t(\eta
_1),(d_{x_2}p_{X'})^t(\eta _2),\eta_1, \eta_2\right ):
\right .\\
&&\left .x_i\in X',\, \eta_i\in T^*_{p_{X'}(x_i)}(X_0),\,
i=1,2.\right \}.
\end{array}
\end{equation}
The fibre diagram is
\begin{equation}\label{eqn:fibdiag1}
  \begin{array}{ccc}
    F_{p_{X'}\times p_{X'}} & \longrightarrow & \Gamma _{p_{X'}\times p_{X'}} \\
     &  & \\
    \downarrow & & \downarrow \\
     &   & \\
    \Sigma '& \longrightarrow & T^*(X'\times X').
  \end{array}
\end{equation}
Here $F_{p_{X'}\times p_{X'}}$ is clearly diffeomorphic to $\Sigma
'$. Thus, the excess is now
\begin{eqnarray*}
  e _{p_{X'}\times p_{X'}}&=&\dim (F_{p_{X'}\times p_{X'}})+
\dim T^*(X'\times X')-\dim \left (\Gamma _{p_{X'}\times p_{X'}}
\right )
  -\dim (\Sigma') \\
  &=& 2g.
\end{eqnarray*}
Given that the fibres of $p_{X'}\times p_{X'}$ have dimension
$2g$, Theorem 9.2 of \textit{loc.cit.} now implies $(p_{X'}\times
p_{X'})_*(\jmath \times \jmath)^*\left (\tilde \Pi _X\right )\in
J^{(1+g)/2}(X_0\times X_0,\Sigma _0)$.

\section{The symplectic structure along the cone.}

Let
\begin{equation}\label{eqn:sigmaalpha}
\Sigma _\alpha =\left \{(x,r\alpha _x)\,:\,x\in X,r>0\right
\}\subseteq T^*(X)\setminus \{0\}\end{equation} be the half-line
bundle generated by the connection 1-form. This is a closed
symplectic cone in $T^*(X)$. Let $\rho:\Sigma _\alpha \rightarrow
X$ and $q=\pi \circ \rho:\Sigma _\alpha \rightarrow M$ be the
projections. Since $\Sigma _\alpha$ is naturally diffeomorphic to
$X\times \mathbb{R}_+$, we have an intrinsic isomorphism of vector
bundles on $\Sigma _\alpha$, $$T(\Sigma _\alpha)\cong \rho
^*(TX)\oplus \mathrm{span}\left \{ \frac{\partial}{\partial
r}\right\},$$ where $\frac{\partial}{\partial r}$ denotes the
generator of the $\mathbb{R}_+$-action on $\Sigma _\alpha$. On the
other hand, the connection 1-form induces a splitting $T(X)=\pi
^*(T(M))\oplus \mathrm{span}\left \{ \frac{\partial}{\partial
\theta}\right\}$, where $\frac{\partial}{\partial \theta}$ denotes
the generator of the $S^1$-action. On the upshot, we have an
isomorphism of vector bundles on $\Sigma _\alpha$:
\begin{equation}
T(\Sigma _\alpha)\cong q^*(T(M))\oplus \mathrm{span} \left \{
\frac{\partial}{\partial r},\frac{\partial}{\partial
\theta}\right\}. \label{eqn:splittingalpha}
\end{equation}
Since $X$ is oriented, it is metalinear and therefore $T^*(X)$ is
metaplectic. Thus $T(T^*(X))$ is a metaplectic vector bundle, and
therefore so is its restriction $\left .T(T^*(X))\right |_{\Sigma
_\alpha}$. We shall now examine the situation in local coordinates
along $\Sigma _\alpha$.

Fix $z=(x,r\alpha _x)\in \Sigma _\alpha$ and set $p=\pi(x)$. Let
$U\subseteq X$ be a coordinate neighbourhood for $X$ near $x$.
Thus $T^*U\subseteq T^*X$ is an open subset, and there is an
obvious symplectomorphism $T^*U\,\cong U\times \mathbb{R}_{2n+1}$,
where the latter is viewed as an open subset of
$\mathbb{R}^{2n+1}\times \mathbb{R}_{2n+1}$, with its natural
symplectic structure (here
$\mathbb{R}_{2n+1}=(\mathbb{R}^{2n+1})^*$). We shall identify
$\mathbb{R}_{2n+1}$ with $\mathbb{R}^{2n+1}$ by means of the
standard scalar product, and view $\left .\alpha \right |_U$ as a
smooth $\mathbb{R}^{2n+1}$-valued form on $U$. Let
$\mathrm{Jac}_y(\alpha)$ be its Jacobian matrix at $y\in U$.

A differential 2-form on $U$ is represented by smooth function
$U\rightarrow A_{2n+1}$, where $A_{2n+1}$ is the space of
antisymmetric matrices of order $2n+1$. Namely, if we denote by
$t_i$ the local coordinates on $U$, the differential form
$\nu=\sum _{i,j}\nu _{ij}dt_i\wedge dt_j$ is represented by the
$A_{2n+1}$-valued function $[\nu _{i,j}]$. The proof of the
following is left to the reader:

\begin{lem}\label{lem:localform}
The matrix valued function $\frac 12\,\left
[\mathrm{Jac}(\alpha)^t-\mathrm{Jac}(\alpha)\right ]$ represents
$\pi ^*(\Omega)$ on $U$. \end{lem}

\noindent In other words, if $y\in U$ and $v=\sum _iv_i\left.
\frac{\partial}{\partial t_i }\right |_y,\, w=\sum _i\left .
w_i\frac{\partial}{\partial t_i }\right |_y\in T_y(X)$ then $$\pi
^*(\Omega)_y(v,w)=\frac 12v^t\left
[\mathrm{Jac}_y(\alpha)^t-\mathrm{Jac}_y(\alpha)\right ]w.$$

\begin{cor} \label{cor:contraction}
Let $V=\sum _iV_i\frac{\partial}{\partial t_i}$ be a vector field
on $U$. Then the contraction $\iota (V)\pi ^*(\Omega)=\pi
^*(\Omega)(V,\cdot)$ is the 1-form represented on $U$ by the
vector valued function $\frac 12\left
[\mathrm{Jac}_y(\alpha)-\mathrm{Jac}_y(\alpha)^t\right
]\,V$.\end{cor}

\begin{lem}\label{lem:splittingalpha}
The pull-back $q^*(TM)$ has a natural metaplectic structure.
\end{lem}

\noindent By (\ref{eqn:splittingalpha}), Lemma
\ref{lem:splittingalpha} implies that $\Sigma$ has a natural
metaplectic structure.

\noindent \textit{Proof.} With some abuse of language, in the
symplectic coordinate chart $U\times \mathbb{R}^{2n+1}$ the
tangent space $T_z (\Sigma _\alpha )$ is the symplectic subspace
of $ \mathbb{R}^{2n+1}\times \mathbb{R}^{2n+1}$ given by
\begin{equation}
  \begin{array}{ccl}
    T_z (\Sigma _\alpha)& =& \left \{\begin{pmatrix}
  v \\
  r\mathrm{Jac}_x(\alpha)v
\end{pmatrix}: v\in \mathbb{R}^{2n+1}\right \}\oplus
\mathrm{span}\left \{\begin{pmatrix}
  0 \\
  \alpha _x
\end{pmatrix}\right \}\\
     & =&\left \{\begin{pmatrix}
  v \\
  r\mathrm{Jac}_x(\alpha)v
\end{pmatrix}: v\in \mathbb{R}^{2n+1}, \alpha _x(v)=0\right \}\\
&&\oplus \,\mathrm{span}\left \{\begin{pmatrix}
  \kappa \\
  r\mathrm{Jac}_x(\alpha)\kappa
\end{pmatrix} ,\begin{pmatrix}
  0 \\
  \alpha _x
\end{pmatrix}\right \},
  \end{array}\label{eqn:splittingalphalocal}
\end{equation}
where $\kappa$ represents $\frac{\partial}{\partial \theta}$.
Clearly, (\ref{eqn:splittingalphalocal}) is just
(\ref{eqn:splittingalpha}) in local coordinates. It follows that
the symplectic orthocomplement of $T_z (\Sigma _\alpha)$ in
$T_z(T^*(X))$ is given by
\begin{equation}\label{eqn:symporthosigmaalpha}
  T_z (\Sigma _\alpha)^\perp\, =\,
\left \{\begin{pmatrix}
  v \\
  r\mathrm{Jac}_x^t(\alpha)v
\end{pmatrix}: v\in \ker(\alpha_x)\subseteq \mathbb{R}^{2n+1}\right
\}.\end{equation} For simplicity, let us use
$q^*(TM)_{\mathrm{hom}}$ as a short hand for the symplectic vector
bundle $\left (q^*(TM),2r\,q^*(\Omega) \right )$ ($r>0$ is the
conic coordinate on $\Sigma_\alpha$). Then
(\ref{eqn:symporthosigmaalpha}) shows that the differential $dq$
induces a symplectic isomorphism between the vector subundle $ T
(\Sigma _\alpha)^\perp$ and the pull-back
$q^*(T(M))_{\mathrm{hom}}$ (i.e., $q^*(TM)_{\mathrm{hom}}$ with
the opposite symplectic structure). On the other hand, rescaling
provides a symplectic isomorphism $q^*(TM)_{\mathrm{hom}}\cong
q^*(TM)$ (the latter as a short hand for $(q^*(TM),q^*(\Omega))$).
. On the upshot, we have a symplectic isomorphism
\begin{equation}\label{eqn:tgbundleonsigma}
  \begin{array}{ccl}
   \left .T(T^*X)\right |_{\Sigma _\alpha} & \cong&q^*(TM)_{\mathrm{hom}}\oplus
q^*(TM)_{\mathrm{hom}}^-\oplus \mathrm{span} \left \{
\frac{\partial}{\partial r},\frac{\partial}{\partial
\theta}\right\} \\
     & \cong&q^*(TM)\oplus
q^*(TM)^-\oplus \mathrm{span} \left \{ \frac{\partial}{\partial
r},\frac{\partial}{\partial \theta}\right\}.
  \end{array}
\end{equation}
Let now $\mathrm{Bp}(q^*(TM))=q^*\big (\mathrm{Bp}(TM)\big )$ and
$\mathrm{Bp}(q^*(TM)^-)$ be the principal $\mathrm{Sp}(n)$-bundles
of all symplectic frames in $q^*(T^*(M))$ and $q^*(T^*(M))^-$,
respectively. Let the automorphism $\sigma
:\mathrm{Sp}(n)\rightarrow \mathrm{Sp}(n)$ be defined by
$$
\sigma (U)\,=:\,\begin{pmatrix}
  I_n & 0 \\
  0 & -I_n
\end{pmatrix}\,U\,\begin{pmatrix}
  I_n & 0 \\
  0 & -I_n
\end{pmatrix}\,\,\,\,\,\,\,\,\,(U\in \mathrm{Sp}(n)).$$
Let us redefine the action of $\mathrm{Sp}(n)$ on
$\mathrm{Bp}(q^*(TM)^-)$ by composing with $\sigma$. Then the map
$\tau:\mathrm{Bp}(q^*(TM))\rightarrow \mathrm{Bp}(q^*(TM)^-)$
given by
$$\tau:(e_1,\ldots,e_n,f_1,\ldots,f_n)\mapsto
(-e_1,\ldots,-e_n,f_1,\ldots,f_n)$$ is an
$\mathrm{Sp}(n)$-equivariant diffeomorphism. We then have an
$\mathrm{Sp}(n)$-equivariant embedding
$\widehat{\tau}:\mathrm{Bp}(q^*(TM))\rightarrow \mathrm{Bp}(\left
.T^*X\right |_{\Sigma _\alpha})$ given by
$$(\mathbf{e},\mathbf{f})\mapsto \left
(\mathbf{e},\mathbf{f},-\mathbf{e},\mathbf{f},\frac{\partial}{\partial
\theta},\frac{\partial}{\partial r}\right ),$$ for
$(\mathbf{e},\mathbf{f})=(e_1,\ldots,e_n,f_1,\ldots,f_n)\in
\mathrm{Bp}(q^*(TM))$. The inverse image of
$\widehat{\tau}:\mathrm{Bp}(q^*(TM))$ in the metaplectic cover of
$\mathrm{Bp}(\left .T^*X\right |_{\Sigma _\alpha})$ is the
asserted metaplectic structure of $q^*(TM)$.

\bigskip

Let us now fix $\zeta =(x,x,r\alpha _x,-r\alpha _x)\in \Sigma$ and
examine the symplectic structure of $T^*(X\times X)$ near $\zeta$.
We shall also identify $T^*(X\times X)$ with $T^*(X)\times
T^*(X)$. In the symplectic product coordinate chart $(U\times
\mathbb{R}^{2n+1})\times (U\times \mathbb{R}^{2n+1})$ the tangent
space $T_\zeta (\Sigma)$ is the isotropic subspace of $\left (
\mathbb{R}^{2n+1}\times \mathbb{R}^{2n+1}\right )\times \left (
\mathbb{R}^{2n+1}\times \mathbb{R}^{2n+1}\right )$ given by
$$T_\zeta (\Sigma)=\left \{\begin{pmatrix}
  v \\
  r\mathrm{Jac}_x(\alpha)v \\
  v \\-rJ_x(\alpha)v
\end{pmatrix}: v\in \mathbb{R}^{2n+1}\right \}\oplus
\mathrm{span}\left \{\begin{pmatrix}
  0 \\
  \alpha _x\\
  0 \\
  -\alpha_x
\end{pmatrix}\right \}.$$
The symplectic annihilator of $T_\zeta (\Sigma)$ is then
\begin{equation}\label{eqn:symplannihilator}
  \begin{array}{ccc}
    T_\zeta (\Sigma)^\perp & =&\left \{\begin{pmatrix}
  v \\
  r\mathrm{Jac}_x(\alpha)^tv \\
  w \\-r\mathrm{Jac}_x(\alpha)^tw
\end{pmatrix}: \alpha _x(v)=\alpha _x(w)=0\right \}
   + T_\zeta (\Sigma).
  \end{array}
\end{equation}

\begin{lem} The sum of vector spaces on the left hand side of
(\ref{eqn:symplannihilator}) is
direct. \label{lem:directsum}
\end{lem}

\noindent \textit{Proof.} Let $a,b,c \in \mathbb{R}$, $v'\in
\mathbb{R}^{2n+1}$ and $v,w\in \ker (\alpha_x)\subseteq
\mathbb{R}^{2n+1}$ be such that
$$
a \,\begin{pmatrix}
  v' \\
  r\mathrm{Jac}_x(\alpha)v' \\
  v' \\-r\mathrm{Jac}_x(\alpha)v'
\end{pmatrix}\,+b\,\begin{pmatrix}
  v \\
  r\mathrm{Jac}_x(\alpha)^tv \\
  w \\-r\mathrm{Jac}_x(\alpha)^tw
\end{pmatrix}\,+c\,\begin{pmatrix}
  0 \\
  \alpha _x\\
  0 \\
  -\alpha_x
\end{pmatrix}=0.$$
If $a$ or $b$ vanish then so do all the other coefficients. We may
otherwise absorb them into $v$ and $w$ sot as to assume $a=-b=1$;
then $v'=v=w\in \ker (\alpha_x)$. Thus we are reduced to the
equality
$$\left (\mathrm{Jac}_x(\alpha)-\mathrm{Jac}_x(\alpha)^t\right )v
=-\frac{c}{r}\,\alpha _x,$$ for a certain $v\in \ker (\alpha_x)$.
Now $\ker (\alpha _x)$ is the horizontal subspace for the
connection, and the skew matrix $\frac 12\left [
\mathrm{Jac}_x(\alpha)^t-\mathrm{Jac}_x(\alpha)\right ]$
represents the 2-form $d_x\alpha=\pi ^*(\Omega)_x$.

\begin{claim}\label{claim:pullbackcontraction}
Let $\mu $ be a 2-form on $M$ and let $W$ be a vector field on
$M$. Denote by $\tilde W$ the horizontal lift of $W$ to a vector
field on $X$, under the given connection. Let $\iota$ be the
contraction operator between vector fields and differential forms.
Then
$$\iota (\tilde W) \, \pi ^*( \mu)\,=\, \pi ^*\left ( \iota
(W)\,\mu\right ).$$
\end{claim}
\noindent In fact, both 1-forms vanish on vertical tangent
vectors, and they obviously take the same values on horizontal
vectors.

Thus, if $v$ represents (in our coordinate patch) the horizontal
lift of a tangent vector $\xi \in T_p(M)$, $p=\pi (x)$, then
$\frac 12\left
[\mathrm{Jac}_x(\alpha)^t-\mathrm{Jac}_x(\alpha)\right ]v$
represents the pull-back of the 1-form $\iota (\xi)\, \Omega_p$
under the differential $d_p\pi :T_x(X)\, \longrightarrow \,
T_p(M)$. But this may not be a multiple of the connection, unless
it vanishes. Thus, $c=0$ and $\left
(\mathrm{Jac}_x(\alpha)-\mathrm{Jac}_x(\alpha)^t\right )v=0$, that
is, $\iota (\xi)\, \Omega _p=0$. This contradicts the
nondegeneracy of $\Omega$, unless $\xi=0$ and therefore $v=0$.
Lemma \ref{lem:directsum} follows.

\bigskip

Therefore, the vector subspace
$$\left \{\begin{pmatrix}
  v \\
  r\mathrm{Jac}_x(\alpha)^tv \\
  w \\-r\mathrm{Jac}_x(\alpha)^tw
\end{pmatrix}: v,w\in \,\ker (\alpha_x)\right \}\subseteq
T_\zeta\left (T^*(X)\times T^*(X)\right )
$$
is symplectomorphic to the symplectic normal space to $\Sigma$ at
$\zeta$, $N_{\Sigma,\zeta}=\left (T_\zeta \Sigma\right )^\perp /
T_\zeta \Sigma$. We shall denote by $\Omega _{\mathrm{can}}$ the
canonical symplectic structure of the cotangent bundles $T^*(X)$
and $T^*(X\times X)\cong T^*(X)\times T^*(X)$ at $(x,r\alpha _x)$
and at $\zeta$, respectively. In our coordinate chart, this is the
canonical symplectic structure on $\left ( \mathbb{R}^{2n+1}\times
\mathbb{R}^{2n+1}\right )$ and $\left ( \mathbb{R}^{2n+1}\times
\mathbb{R}^{2n+1}\right )\times \left ( \mathbb{R}^{2n+1}\times
\mathbb{R}^{2n+1}\right )$. We have, for $v,w,v',w' \in \ker
(\alpha _x)$ horizontal lifts of $\xi,\eta,\xi',\eta'\in T_pM$:
\begin{equation}
  \begin{array}{lc}
\Omega _{\mathrm{can}}\left (
\begin{pmatrix}
  v \\
  r\mathrm{Jac}_x(\alpha)^tv \\
  w \\-r\mathrm{Jac}_x(\alpha)^tw
\end{pmatrix},\begin{pmatrix}
  v' \\
  r\mathrm{Jac}_x(\alpha)^tv' \\
  w' \\-r\mathrm{Jac}_x(\alpha)^tw'
\end{pmatrix}\right )= &  \\
  \Omega _{\mathrm{can}}\left (
\begin{pmatrix}
  v \\
  r\mathrm{Jac}_x(\alpha)^tv
\end{pmatrix},\begin{pmatrix}
  v' \\
  r\mathrm{Jac}_x(\alpha)^tv'
\end{pmatrix}\right )-\,\Omega _{\mathrm{can}}\left (
\begin{pmatrix}
  w \\
  r\mathrm{Jac}_x(\alpha)^tw
\end{pmatrix},\begin{pmatrix}
  w' \\
  r\mathrm{Jac}_x(\alpha)^tw'
\end{pmatrix}\right )& \\
=r\Omega _p(\xi,\xi')-r\Omega _p(\eta,\eta').
  \end{array}
\end{equation}
Thus, there is a natural symplectic isomorphism
\begin{equation}\label{eqn:sympnormalbundles}
N_{\Sigma,\zeta}\, \cong \, (T_pM,\,2r\Omega _p) \oplus
(T_pM,\,-2r\Omega _p),\end{equation} which extends to an
isomorphism of symplectic vector bundles with the appropriate
symplectic structures. Let $J_M$ be the complex structure of $M$.
Then, when endowed with the compatible complex structure
$(J_M,-J_M)$, $N_{\Sigma}$ is a unitary (hence Riemannian) vector
bundle.

Let us next consider the vector subspace $V_\zeta \subseteq
T_\zeta(T^*(X\times X))\setminus \{0\}$ given by:
$$V_\zeta=:\mathrm{span}\left \{
\begin{pmatrix}
 v \\
 r\,\mathrm{Jac}_x(\alpha)v \\
 w\\
-r\mathrm{Jac}_x(\alpha)w
\end{pmatrix}\,:\,
\alpha _x(v)=\alpha (w)=0\right \}.$$ Then $V_\zeta$ is also
naturally symplectomorphic to $(T_pM,-2r\Omega_p)\oplus
(T_pM,2r\Omega_p)$, $p=\pi (x)$, and thus becomes a unitary vector
space with the complex structure $(-J_p,J_p)$; it contains
$T_\zeta(\Sigma)$ as a Lagrangian subspace. Clearly, it naturally
extends to a vector subbundle of $V\subseteq \left .T(T^*(X)\times
T^*(X)\right )$.

Next, let $\kappa_\theta:U\rightarrow \mathbb{R}^{2n+1}$ represent
the vector field $\frac{\partial}{\partial \theta}$ generating of
the $S^1$-action, and let us set $\kappa =\kappa _\theta (x)$. Let
us then introduce the vector subspace $W_\zeta \subseteq
T_\zeta(T^*(X)\times T^*( X))\setminus \{0\}$ given by
$$W_\zeta =\mathrm{span}\left \{\begin{pmatrix}
 \kappa \\
  r\mathrm{Jac}_x(\alpha)\kappa \\
  0
  \\
  0
\end{pmatrix},\begin{pmatrix}
  0 \\
  r\alpha \\
  0 \\
  0
\end{pmatrix},\begin{pmatrix}
  0 \\
 0 \\
  0 \\
  -r\alpha
\end{pmatrix},\begin{pmatrix}
 0 \\
  0\\
  \kappa
  \\
  -r\mathrm{Jac}_x(\alpha)\kappa
\end{pmatrix}\right \}.
$$
Since $\kappa ^t\cdot \alpha _x=\alpha _x\left ( \left
.\frac{\partial}{\partial \theta}\right |_x\right )=1$, by mapping
this (intrisically defined) basis to the real basis
$(e_1,ie_1,e_2,ie_2)$ of $\mathbb{C}^2$ we see that $W_\zeta$ is
naturally symplectomorphic to $(\mathbb{C}^2,r\omega_0)$, where
$\omega_0=\frac{i}{2}\,\sum _{i=1}^2dz_i\wedge d\overline z_i$ is
the standard symplectic structure. We shall then consider
$W_\zeta$ as the fibre of a unitary vector bundle $W$ on $\Sigma$.
The proof of following is left to the reader:
\begin{lem}\label{lem:tgbundlexperxsigma}
We have the symplectic direct sum decomposition
$$
T_\zeta\left (T^*(X)\times T^*(X)\right )\,\cong\,
N_{\Sigma,\zeta} \oplus V_\zeta\oplus W_\zeta.$$
\end{lem}
We may then take the direct sum of the unitary structures on each
summand to make $ \left .T\left (T^*(X)\times T^*(X)\right )\right
|_\Sigma$ into a unitary vector bundle over $\Sigma$.

We next consider the closed isotropic cone
$$
\Sigma '=\left \{(x,r\alpha '_x,x,-r\alpha '_x)\,:\,x\in
X',r>0\right \}\subseteq T^*(X')\times T^*(X')\setminus \{0\},$$
where $\alpha'=\jmath ^*(\alpha)$. Let us fix $\zeta '=(x,r\alpha
'_x,x,-r\alpha _x')\in \Sigma'$, and choose a coordinate patch
$U'\subseteq X'$ containing $x$. Thus $T^*(U')\cong U'\times
\mathbb{R}^{2n+1-g}$, and $\alpha '$ is represented locally near
$x$ by a smooth function $U'\rightarrow \mathbb{R}^{2n+1-g}$. Let
$\mathrm{Jac}'_x(\alpha ')$ be its Jacobian matrix at $x$. Then in
local coordinates
$$
T_{\zeta'} (\Sigma ')=\left \{
\left(%
\begin{array}{c}
  v \\
  r\mathrm{Jac}'_x(\alpha')v \\
  v \\
  -r\mathrm{Jac}'_x(\alpha')v \\
\end{array}%
\right)\,:\,v\in \mathbb{R}^{2n+1-g}\right \}\oplus
\mathrm{span}\left \{\left(%
\begin{array}{c}
  0 \\
  \alpha '_x\\
  0 \\
  -\alpha '_x\\
\end{array}%
\right) \right \}.$$ The symplectic annihilator is then
$T_{\zeta'} (\Sigma ')^\perp=R_{\zeta'}+T_{\zeta'} (\Sigma ')$,
where
\begin{equation}\label{eqn:sympann'}
\begin{array}{ccc}
   R_{\zeta'} & =&\left \{
\left(%
\begin{array}{c}
  v \\
  r\mathrm{Jac}'_x(\alpha')^tv \\
  w \\
  -r\mathrm{Jac}'_x(\alpha')^tw \\
\end{array}%
\right)\,:\,v,w\in \ker (\alpha '_x)\right \} \\
     & & \oplus\,\mathrm{span}\left \{
\left(%
\begin{array}{c}
   0\\
  g_{\xi,x} \\
  0 \\
  -g_{\xi,x} \\
\end{array}%
\right):\xi \in \frak{g}\right \}.
  \end{array}
\end{equation}
Here, for every $\xi \in \frak{g}$, $g_{\xi,x}$ denotes the linear
functional on $T_x(X')$ given by the Riemannian scalar product
with $\xi ^\sharp(x)$, where $\xi ^\sharp$ is the vector field on
$X'$ generated  by $\xi$. Because of the degeneracy of the
restricted form, the two summands for $T_{\zeta'} (\Sigma
')^\perp$ have a $g$-dimensional intersection:

\begin{lem}\label{lem:intersR}
$$R_{\zeta'} \cap T_{\zeta'} (\Sigma ')= \left \{
\left(%
\begin{array}{c}
   \xi ^\sharp (x)\\
  r\mathrm{Jac}'_x(\alpha ')\xi ^\sharp (x) \\
  \xi ^\sharp (x) \\
  -r\mathrm{Jac}'_x(\alpha ')\xi ^\sharp (x) \\
\end{array}%
\right):\xi \in \frak{g}\right \}.$$
\end{lem}
\noindent \textit{Proof.} Notice to begin with that the action of
$G$ on $X$ is horizontal on $X'$, that is, $\xi ^\sharp (x)\in
\ker (\alpha _x')$ if $\xi \in \frak{g}$ and $x\in X'$. Let
$\Omega '=\iota ^*\big (\Omega\big )$, where $\iota
:M'\hookrightarrow M$ is the inclusion. Then $\Omega'$ is a
degenerate closed 2-form, whose kernel at any $p\in M'$ is the
tangent space $\frak{g}\cdot p\subseteq T_p(M)$ to the orbit
through $p$. On the other hand, $\pi ^{\prime *}\left
(\Omega'\right )=d\alpha '$, where $\pi ':X'\rightarrow M'$ is the
projection. Since the space of all $\xi ^\sharp (x)$'s is the
horizontal lift of $\frak{g}\cdot p$ at $x$, we deduce that in our
local coordinates we have $\xi ^\sharp(x)\in \ker\left
(\mathrm{Jac}_x'(\alpha')-\mathrm{Jac}_x'(\alpha')^t\right )$ for
all $\xi \in \frak{g}$. On the upshot, $\ker (\alpha '_x)\cap
\ker\left
(\mathrm{Jac}_x'(\alpha')-\mathrm{Jac}_x'(\alpha')^t\right )$ is
precisely the space of all $\xi ^\sharp (x)$'s, $\xi \in
\frak{g}$. Suppose now that we have an equality:
$$
\left(%
\begin{array}{c}
  v' \\
  r\mathrm{Jac}'_x(\alpha')v' \\
  v' \\
  -r\mathrm{Jac}'_x(\alpha')v' \\
\end{array}%
\right)+b
\left(%
\begin{array}{c}
   0\\
  \alpha  '_x\\
  0 \\
  -\alpha '_x \\
\end{array}%
\right)\,=\,
\left(%
\begin{array}{c}
  v \\
  r\mathrm{Jac}'_x(\alpha')^tv \\
  w \\
  -r\mathrm{Jac}'_x(\alpha')^tw \\
\end{array}%
\right)+
\left(%
\begin{array}{c}
   0\\
  g_{\xi,x} \\
  0 \\
  -g_{\xi,x} \\
\end{array}%
\right),$$ with $v,w\in \ker (\alpha _x')$ and $\xi \in \frak{g}$.
Then $v'=v=w\in \ker (\alpha '_x)$ and so:
$$r(\mathrm{Jac}'_x(\alpha')-\mathrm{Jac}'_x(\alpha')^t)v+ b\alpha
'_x- g_{\xi,x}=0.$$  The left hand side is a cotangent vector to
$X'$ at $x$. By pairing this first with the generator at $x$ of
the $S^1$-action on $X'$, and then with the $\eta ^\sharp (x)$'s,
$\eta \in \frak{g}$, we obtain $b=0$ and $\xi=0$. Thus, $v=\xi
^\sharp (x)$ for some $\xi \in \frak{g}$.

\bigskip

Let $H(M'/M_0)\subseteq TM'$ be the Riemannian orthocomplement of
the vertical tangent bundle of $p_{M'}:M'\rightarrow M_0$. Then
$H(M'/M_0)$ is a connection for the principal $G$-bundle $p_{M'}$.
When $M'$ is endowed with the 2-form $\Omega'=\iota
^*(\Omega)=p_{M'} ^*(\Omega _0)$, $H(M'/M_0)$ is a symplectic
vector subbundle of $TM'$, symplectomorphic to the pull-back
$p_{M'}^*\left (TM_0\right )$. Let $\Omega _{H(M'/M_0)}$ be its
symplectic structure. Now $R_\zeta$ in (\ref{eqn:sympann'}) may be
decomposed as
\begin{equation}\label{eqn:sympann''}
\begin{array}{ccl}
   R_{\zeta'} & =&\left \{
\left(%
\begin{array}{c}
  v \\
  r\mathrm{Jac}'_x(\alpha')^tv \\
  w \\
  -r\mathrm{Jac}'_x(\alpha')^tw \\
\end{array}%
\right)\,:\,v,w\in \ker (\alpha '_x),\, d_x\pi (v),d_x\pi (w)\in H_p\right \} \\
& &\oplus\left \{
\left(%
\begin{array}{c}
  \xi ^\sharp (x) \\
  r\mathrm{Jac}'_x(\alpha')^t\xi ^\sharp (x) \\
  -\xi ^\sharp (x) \\
  r\mathrm{Jac}'_x(\alpha')^t\xi ^\sharp (x) \\
\end{array}%
\right)\,:\,\xi\in \frak{g}\right \}\oplus\,\mathrm{span}\left \{
\left(%
\begin{array}{c}
   0\\
  g_{\xi,x} \\
  0 \\
  -g_{\xi,x} \\
\end{array}%
\right):\xi \in \frak{g}\right \}\\
& &\oplus\left \{
\left(%
\begin{array}{c}
  \xi ^\sharp (x) \\
  r\mathrm{Jac}'_x(\alpha')^t\xi ^\sharp (x) \\
  \xi ^\sharp (x) \\
  -r\mathrm{Jac}'_x(\alpha')^t\xi ^\sharp (x) \\
\end{array}%
\right)\,:\,\xi\in \frak{g}\right \}.
  \end{array}
\end{equation}

By Lemma \ref{lem:intersR}, the first two summands add up
symplectomorphically to the symplectic normal bundle to $\Sigma'$
in $T^*(X')\times T^*(X')$. Thus,

\begin{lem}
Let $N_{\Sigma '}'$ be the symplectic normal bundle of $\Sigma
'\subseteq T^*(X')\times T^*(X')$. Let $\zeta'
=(x,r\alpha'_{x'},x,-r\alpha '_{x'})$, $p=\pi (x)\in M$,
$\overline p=q(p)\in M_0$. Then we have a natural isomorphism of
symplectic vector spaces:
\begin{equation*}
  \begin{array}{ccl}
   N'_{\Sigma',\zeta'}& =&(H_p(M'/M_0),2r\,\Omega _{H(M'/M_0),p})\oplus (H_p(M'/M_0),-2r\,\Omega _{H(M'/M_0),p})\oplus
\frak{g}_{\mathbb{C},\overline p} \\
     & \cong&(T_{\overline p}(M_0),r\Omega _{0,\overline p})\oplus (T_{\overline p}(M_0),
-r\Omega _{0,\overline p})\oplus \frak{g}_{\mathbb{C}, p}.
  \end{array}
\end{equation*}
This naturally extends to an isomorphism of unitary vector bundles
over $\Sigma'$.
\end{lem}

Here, $$\frak{g}_{\mathbb{C}, p}=\mathrm{span}\left
\{\left(%
\begin{array}{c}
   \xi ^\sharp (x)\\
  r\mathrm{Jac}'_x(\alpha') \xi ^\sharp (x)\\
  -\xi ^\sharp (x) \\
  r\mathrm{Jac}'_x(\alpha')\xi ^\sharp (x)\\
\end{array}%
\right),\left(%
\begin{array}{c}
   0\\
  g_{\xi,x} \\
  0 \\
  -g_{\xi,x} \\
\end{array}%
\right):\xi \in \frak{g}\right \}.$$ On $T_p(M)\supseteq T_p(M')$,
we have in obvious notation $g_\xi =\Omega _p(\cdot,J_p\xi ^\sharp
(p))$. Thus, $ \frak{g}_{\mathbb{C}, p}$ may be naturally
identified with the complexified Lie algebra of $G$, that is, the
Lie algebra of the complexified group $\tilde G$. It is endowed
with the unitary structure induced by its infinitesimal action on
$M$ at $p$, which makes it into a complex subspace of $T_p(M)$:
$\frak{g}_{\mathbb{C}, p}\cong \frak{g}_{\mathbb{C}}\cdot  p=T_{
p}(\tilde G\cdot p)$. Let $\kappa '_\theta$ denote the generator
of the $S^1$-action on $X'$, set $\kappa '=\kappa'_\theta(x)$ and
let us now introduce the vector bundles over $\Sigma'$ by setting,
in local coordinates,
\begin{equation*}
  \begin{array}{ccl}
    V'_{\zeta'} & =&\left \{\begin{pmatrix}
     v \\
      r\mathrm{Jac}_x(\alpha ')v \\
      w \\
      -r\mathrm{Jac}_x(\alpha ')w
    \end{pmatrix}:v,w\in \ker (\alpha '_x),\, d_x\pi(v), d_p\pi (w)\in H_p\right
\},
  \end{array}
\end{equation*}
\begin{equation*}
  \begin{array}{ccl}
    V''_{\zeta'} & =&\left \{\begin{pmatrix}
     \xi ^\sharp (x) \\
      r\mathrm{Jac}_x(\alpha ')\xi ^\sharp(x) \\
      \xi^\sharp (x) \\
      -r\mathrm{Jac}_x(\alpha ')\xi ^\sharp
    \end{pmatrix},\,\begin{pmatrix}
     0 \\
      g_\xi \\
      0 \\
      g_\xi
    \end{pmatrix}:\xi \in \frak{g}\right
\},
  \end{array}
\end{equation*}
$$W'_{\zeta'} =\mathrm{span}\left \{\begin{pmatrix}
 \kappa '\\
  r\mathrm{Jac}'_x(\alpha')\kappa '\\
  0
  \\
  0
\end{pmatrix},\begin{pmatrix}
  0 \\
  r\alpha '_x\\
  0 \\
  0
\end{pmatrix},\begin{pmatrix}
  0 \\
 0 \\
  0 \\
  -r\alpha_x'
\end{pmatrix},\begin{pmatrix}
 0 \\
  0\\
  \kappa'
  \\
  -r\mathrm{Jac}'_x(\alpha')\kappa'
\end{pmatrix}\right \}.
$$
Just as before, these are naturally unitary vector bundles, and we
have:
\begin{lem} \label{lem:tgbundlexprxxpr}There is
a symplectic direct sum decomposition
$$
\left .T(T^*(X)\times T^*(X))\right |_{\Sigma'}\cong N'_{\Sigma'}
\oplus V'\oplus V''\oplus W'.$$
\end{lem}
By taking the direct sum of the unitary structure of each summand,
this makes $ \left .T(T^*(X)\times T^*(X))\right |_{\Sigma'}$ into
a unitary vector bundle.

\section{The symbol of $(\jmath \times \jmath)^*\left (\tilde \Pi _{X}\right )$.}

If $W$ is a manifold, the symbol of a Fourier-Hermite distribution
$u\in J^k(W,\Xi)$ associated to a closed isotropic cone $\Xi
\subseteq T^*W\setminus \{0\}$ is a smooth section of the
symplectic spinor bundle, $\mathrm{Spin}(\Xi)$, homogeneous of
degree $k$ with respect to the conic structure of $\Xi$. The space
of homogeneous sections of degree $k$ of $\mathrm{Spin}(\Xi)$ will
be denoted by $S^k(\Xi)$.

In particular, the symbol of the Szeg\"{o} kernel $\tilde \Pi
_X\in J^{1/2}(X\times X,\Sigma)$ is an element of
$S^{1/2}(\Sigma)$. The spinor bundle of $\Sigma \subseteq
T^*(X\times X)\subseteq \{0\}$ is $\mathrm{Spin}(\Sigma)=\bigwedge
^{1/2}\Sigma \otimes \mathrm{S}(N_\Sigma)$, where $N_\Sigma$ is
the symplectic normal bundle of $\Sigma$. Now $\Sigma$ is
obviously diffeomorphic to $\Sigma _\alpha$ in
(\ref{eqn:sigmaalpha}), and by the symplecticity of $\Omega$ the
latter is a symplectic submanifold of $T^*X$.

\begin{defn}
By the above, $\Sigma$ carries a built-in symplectic structure
homogeneous of degree one. We shall denote by $\mathrm{vol}_\Sigma
^{1/2}$ the nowhere vanishing half-form of degree $1/2$ obtained
from the latter by the appropriate homogenization with respect to
the $r$ coordinate.\end{defn}

On the other hand, in view of (\ref{eqn:sympnormalbundles}), the
symplectic normal bundle of $\Sigma$ is (after an appropriate
rescaling) naturally isomorphic to $q^*(TM)\oplus q^*(TM)^-$.
Thus, given the complex structure $J_M$ of the base manifold $M$,
by Corollary \ref{cor:sympdirectsumpm} the bundle
$\mathcal{S}(N_\Sigma)\in \mathrm{End}_{\mathcal{HS}}\left
(\mathcal{S}(q^*(TM)\right )$ has a built-in nowhere vanishing
section $\sigma _{J_M}$.

 After \cite{bg}, the symbol of the Szeg\"{o} kernel is the tensor product
\begin{equation}\label{eqn:szegosymbol}
\sigma \left (\tilde \Pi _X\right )\, =\,\mathrm{vol}_\Sigma
^{1/2}\otimes \sigma _{J_M}.\end{equation}

As an intermediate step towards computing the symbol of $\tilde
\Pi _{X|X_0}$, we shall now compute the symbol of the restriction
$(\jmath \times \jmath)^*\left (\tilde \Pi _{X}\right )$.

If $x\in X'\subseteq X$ and $\eta \in T^*_xX$ is a cotangent
vector to $X$ at $x$, we shall use the notation $\eta '=\left
(d_x\jmath\right )^*(\eta)\in T^*_x(X')$ for the restriction of
$\eta$ to $X'$. The relevant canonical relation in $T^*(X'\times
X')\times T^*(X\times X)$ is thus
\begin{equation*}
\begin{array}{ccl}
  \Gamma _{\jmath\times \jmath} & =&\left \{\big ((x_1,\eta'_1,x_2,\eta_2'),
(x_1,\eta_1,x_2,\eta_2)\big):x_i\in X', \, \eta_i\in
T_{x_i}^*(X)\right \},
\end{array}
\end{equation*}
and the fibre product $F_{\jmath\times \jmath}$ of $\Sigma$ and
$\Gamma _{\jmath\times \jmath}$ maps diffeomorphically to
$\Sigma'$ under the projection $p_{F_{\jmath\times
\jmath}}:F_{\jmath\times \jmath}\rightarrow \Gamma _{\jmath\times
\jmath}\circ \Sigma=\Sigma'$. Let $q:T^*(X'\times X')\times
T^*(X\times X)\rightarrow T^*(X\times X)$ be the projection onto
the second factor. Fix $\zeta =(x,r\alpha,x,-r\alpha)\in \Sigma$;
the differential of $q$ at $(\zeta',\zeta)$,
$$q_\zeta=:d_{(\zeta',\zeta)}q: T_{\zeta'}\left (T^*(X'\times
X')\right )\times T_\zeta\left (T^*(X\times X)\right ) \rightarrow
T_\zeta\left (T^*(X\times X)\right ),$$ is simply projection onto
the second factor. Then we have (in local coordinates)
\begin{equation*}
\begin{array}{ccl}
  q_\zeta\left (T _{(\zeta',\zeta)}(\Gamma _{\jmath\times \jmath})\right ) & =&\left \{
\begin{pmatrix}
  v \\
  \phi \\
 w \\
  \eta
\end{pmatrix}\,:\,v,w \in T_x(X'),\phi,\eta \in \mathbb{R}^{2n+1}\right
\}.
\end{array}
\end{equation*}
Given that $X'=(\Phi \circ \pi)^{-1}(0)$, its tangent space
$T_x(X')$ is defined by the vanishing of all the differentials
$d_xH_\xi$, $\xi \in \frak{g}$, where $H_\xi=<\xi,\Phi>$ is the
Hamiltonian function associated to $\xi$. Its symplectic
annihilator in $T_\zeta\big (T^*(X)\times T^*(X)\big )$ is then
\begin{equation*}
\begin{array}{ccl}
  q_\zeta\left (T _{(\zeta',\zeta)}(\Gamma _{\jmath\times \jmath})\right )^{\perp} & =&\left \{
\begin{pmatrix}
  0 \\
 d_xH_\xi \\
 0 \\
 d_xH_\eta
\end{pmatrix}\,:\,\xi, \eta \in \frak{g}\right
\}.
\end{array}
\end{equation*}
It follows that $U_{0\zeta}=:q_\zeta\left (T
_{(\zeta',\zeta)}(\Gamma _{\jmath\times \jmath})\right
)^{\perp}\cap T_\zeta(\Sigma)=\{0\}$, while
$$U_{1\zeta}=:q_\zeta\left (T _{(\zeta',\zeta)}(\Gamma
_{\jmath\times \jmath})\right )^{\perp}\cap
T_\zeta(\Sigma)^{\perp}=\mathrm{span}\left \{\begin{pmatrix}
  0 \\
 d_xH_\xi \\
 0 \\
 -d_xH_\xi
\end{pmatrix}\,:\,\xi \in \frak{g}\right
\}.
$$
Thus, in view of the equality holding in local coordinates on
$U\subseteq X$ (see Claim \ref{claim:pullbackcontraction})
$$dH_\xi=\iota
(X_\xi)\Omega\,=\, \frac 12\,
(\mathrm{Jac}(\alpha)-\mathrm{Jac}(\alpha)^t)X_\xi\,\,\,\,\,(\xi
\in \frak{g}),$$ the vector space $U_{1\zeta}\cong
U_\zeta\subseteq N_{\Sigma,\zeta}\cong (T_pM,2r\Omega_p)\oplus
(T_pM,-2r\Omega_p)$ may be identified with the isotropic subspace
\begin{equation}\label{eqn:Udiagonal}
U_p\cong \left \{(X_\xi(p),X_\xi(p)):\xi \in \frak{g}\right
\},\end{equation} where $X_\xi$ us the vector field generated by
$\xi \in \frak{g}$. As above, we shall make
$(T_pM,2r\Omega_p)\oplus (T_pM,-2r\Omega_p)$ into a unitary vector
space, by endowing it with the compatible complex structure
$(J_p,-J_p)$. Furthermore, let $\frak{g}^M$ denote the trivial
vector subbundle of $\left .TM\right |_{M'}$ on $M'$ with fibre
$\frak{g}$ generated by the infinitesimal action of $\frak{g}$ on
$M$. Then $\frak{g}^M$ is an oriented isotropic subbundle of
$\left .TM\right |_{M'}$. We may write (\ref{eqn:Udiagonal}) as
\begin{equation}\label{eqn:diagonalU}
  U=\mathrm{diag}\left (\frak{g}^M\right )\subseteq (T_pM,2r\Omega_p)\oplus
(T_pM,-2r\Omega_p).
\end{equation}
Therefore, Corollary \ref{cor:imageEpmsurjmorph} may be applied to
the symplectic normal bundle $N_U$ of $U$ in
$(T_pM,2r\Omega_p)\oplus (T_pM,-2r\Omega_p)$. To this end, let us
adopt the notation introduced in the discussion preceding
(\ref{eqn:sympann''}) and in Corollary \ref{cor:imageEpmsurjmorph}
with $L=\frak{g}^M$, and recall after  \cite{gs-gq} that the
symplectic and Riemannian annihilators of $\frak{g}^M$ within
$\left .TM\right |_{M'}$ are given by $\left (\frak{g}^M\right
)^\perp=\frak{g}^M\oplus H(M'/M_0)$ and $\left (\frak{g}^M\right
)^0=J_M(\frak{g}^M)\oplus H(M'/M_0)$, respectively.

Let us then define vector bundles $\frak{g}^M_r$ and
$\frak{g}^M_i$ on $M'$ by setting
$\frak{g}^M_r(p)=:\{(X_\xi(p),-X_\xi(p)):\xi \in \frak{g}\}$ and
$\frak{g}^M_i(p)=:\{(J_pX_\xi(p),J_pX_\xi(p)):\xi \in \frak{g}\}$
($p\in M'$), and then set
$(\frak{g}^M)_{\mathbb{C}}=\frak{g}^M_r\oplus \frak{g}^M_i$.

By Corollary \ref{cor:imageEpmsurjmorph} iii), we conclude
\begin{cor}\label{cor:normalbundleU}
the symplectic annihilator of $U$ in $N_\Sigma$ is
\begin{equation}\label{eqn:normalU}
\begin{array}{ccl}
 N_U&\cong & (\frak{g}^M)_{\mathbb{C}}\\
&& \oplus (q^*H(M'/M_0),2r\Omega
 _{H(M/M')})\oplus (q^*H(M'/M_0),-2r\Omega _{H(M/M')}).
\end{array}\end{equation}
\end{cor}

\begin{lem} $N_U$ is a metaplectic vector bundle.
\end{lem}

\noindent \textit{Proof.} Let us consider the nested fibre
diagrams
\begin{equation}
  \begin{array}{cccccccc}
   &\Sigma'&\stackrel{\rho'}{\longrightarrow}  & X'& \stackrel{\pi_{X'}}{\longrightarrow}&M' &\\
   &&  & & & &\\
p_{\Sigma'}&\downarrow&p_{X'}&\downarrow& &\downarrow& p_{M'}\\
&&&&\\
  &\Sigma _0& \stackrel{\rho_0}{\longrightarrow}    & X_0& \stackrel{\pi_{X_0}}{\longrightarrow}&M_0.&
  \end{array}
\end{equation} The horizontal arrows are principal $S^1$-bundles,
and the vertical ones are principal $G$-bundles. Set $q'=\pi
_{X'}\circ \rho ':\Sigma'\rightarrow M'$, $q_0=\pi_{X_0}\circ \rho
_0:\Sigma _0\rightarrow M_0$. It suffices to show that $q^{\prime
* }(H(M'/M_0))$ is a metaplectic vector bundle, and thus by commutativity that
$q_0^*(TM_0)$ is metaplectic. This follows from Lemma
\ref{lem:splittingalpha} applied to $q_0$.

Equivalently, $\left .q^*(H(M'/M_0))\right |_{\Sigma'}$ is the
symplectic normal bundle of the oriented isotropic subbundle
$\left .q^*\left (\frak{g}^M\right )\right |_{\Sigma'}$ in the
metaplectic vector bundle $\left .q^*(TM)\right |_{\Sigma'}$, and
the fact that is is metaplectic then follows from Proposition
\ref{prop:metaplbiject}.

\bigskip

In view of Lemma \ref{lem:usefulfacts} and Corollary
\ref{cor:sympdirectsumpm}, we have
$$\mathcal{S}(N_U)\cong \mathcal{S}\left (
\frak{g}_{\mathbb{C}}^M\right )\otimes
\mathrm{End}_{\mathcal{HS}}\left (\mathcal{S}(q^{\prime *
}H(M'/M_0))\right ).$$ Given that $(\frak{g}^M)_{\mathbb{C}}$ is a
complex vector subbundle of $TM$, the complex structure $J_M$ of
$M$ singles out the section $\sigma _{\frak{g}^M,J_M}$ of
$\mathcal{S}((\frak{g}^M)_{\mathbb{C}}))$. Since $H(M'/M)\cong
p_{M'}^*(TM_0)$, the complex structure $J_0$ of $M_0$ singles out
the section $\sigma _{J_0}$ of $\mathrm{End}_{\mathcal{HS}}\left
(\mathcal{S}(q^{\prime * }H(M'/M_0))\right )$. Let $\sigma _{J_M}$
be as in (\ref{eqn:szegosymbol}) and let $\Phi _U:
\mathcal{S}(N_\Sigma)\rightarrow \bigwedge ^{-1/2}(U)\otimes
\mathcal{S}(N_U)$ be the morphism of vector bundles from Corollary
\ref{cor:surjmorphi}. Then, in view of Corollary
\ref{cor:imageEpmsurjmorph}, iv), we have
\begin{equation}\label{eqn:fromnsigmatonU}
  \Phi (\sigma _{J_M})=\mathrm{vol}^{-1/2}(U)\otimes \sigma
_{\frak{g}^M,J_M}\otimes \sigma _{J_0},
\end{equation}
where $\mathrm{vol}^{-1/2}(U)$ is the $-\frac 12$-form on $U$
taking value one on oriented orthonormal frames of $U$.

Let now $\mathrm{vol}_{\Gamma_{\tilde \jmath \times \tilde
\jmath}} ^{1/2}$ be the half-form on $\Gamma_{\tilde \jmath \times
\tilde \jmath}$ associated to the morphism of metalinear manifolds
$\tilde \jmath \times \tilde \jmath$ as in
(\ref{eqn:halfformmeta}), and let $\mathrm{vol}_\Sigma ^{1/2}$ be
as in (\ref{eqn:szegosymbol}). In the following definition, recall
that the complex structure $J_M$ maps the vector subbundle
$\frak{g}^M$ of $\left .TM\right |_{M'}$ isometrically onto the
vector subbundle $\frak{g}^M$.

\begin{defn} \label{def:volumedisigma'}
Let $\varepsilon$ be the half-form on $\frak{g}^M$ taking value
one on oriented orthonormal frames. Given the exact sequence
$$0\longrightarrow T(\Sigma') \longrightarrow
\left .T(\Sigma)\right |_{\Sigma'}\longrightarrow
J_M(\frak{g}^M)\longrightarrow 0,$$ let
$\mathrm{vol}_{\Sigma'}^{1/2}$ be the half-form on $\Sigma'$,
homogeneous of degree $\frac 12$, obtained by dividing
$\mathrm{vol}_\Sigma ^{1/2}$ by $\varepsilon$.
\end{defn}

\begin{lem}\label{lem:imageunderlbmap}
Suppose $x\in X$ and let $\zeta =(x,\alpha _x,x,-\alpha _x)$
(thus, $r=1$). Then the image at $\zeta$ of
$\mathrm{vol}_{\Gamma_{\tilde \jmath \times \jmath},\zeta}
^{1/2}\otimes \mathrm{vol}_{\Sigma,\zeta }^{1/2}\otimes
\mathrm{vol}_{U,\zeta} ^{-1/2}$ in $\bigwedge
^{1/2}(\Sigma')_\zeta$ under the line bundle isomorphism
(\ref{eqn:isomolnbdl}) is
$2^{g}\,\mathrm{vol}_{\Sigma',\zeta}^{1/2}$.
\end{lem}

The factor $2^g$ is related to the relative dimension of the
embedding $\tilde j\times \tilde j$, which is $2g$. A similar
factor of $2^{-g}$ will appear in the following section, due to
the fact that the relative dimension of $p_{X'}\times p_{X'}$ is
$-2g$. The two factors will thus cancel out in the final result.

\noindent \textit{Proof.} We shall henceforth omit the base point
$\zeta$. Let us fix an oriented orthonormal (real) basis
$\{v_j\}_{1\le j\le 2n}$ of $(T_pM,2r\Omega_p)$, where $p=\pi
(x)$, of the following form. First, let us choose an oriented
orthonormal basis $\{v_j\}$ of $H_p(M'/M_0)$ for $1\le j\le
2(n-g)$. Next, let us set $v_{2(n-g)+j}=\xi _j^\sharp(p)$, $1\le
j\le g$, where $\{\xi _j^\sharp(p)\}$ is any orthonormal frame of
$\frak{g}^M(p)$ (the vectors $\xi _j\in \frak{g}$ depend on $p$),
and finally let us set $v_{2n-g+j}=J_{M,p}\xi _j^\sharp(p)$, $1\le
j\le g$. In particular, $\{v_j\}_{1\le j\le 2n-g}$ is an oriented
orthonormal basis of $T_pM'$. To simplify our expressions, we
shall in the following denote by $\overrightarrow{\mathbf{v}}$ the
sequence
$$\left \{ \frac{1}{\sqrt{2}} v_1,\ldots ,\frac{1}{\sqrt{2}}
v_{2(n-g)} \right \},$$ by $\overrightarrow{\mathbf{\xi ^\sharp}}$
the sequence
$$\left \{ \frac{1}{\sqrt{2}} \xi _1^\sharp(p),\ldots
,\frac{1}{\sqrt{2}} \xi _g^\sharp(p) \right \},$$ and by
$\overrightarrow{\mathbf{J\xi ^\sharp}}$ the sequence
$$\left \{ \frac{1}{\sqrt{2}} J_{M,p}X_{\xi_1}(p),\ldots
,\frac{1}{\sqrt{2}} J_{M,p}\xi _g^\sharp(p) \right \}.$$ With our
usual slight abuse of language, we shall identify $\Sigma '$ with
the subset of $\Sigma$ lying over $\mathrm{diag}(X')$. As in
section 5, we shall denote by $\kappa$ the generator of the
$S^1$-action on $X$ (expressed in local coordinates). As an
oriented orthonormal basis of $T_\zeta(\Sigma')$, we shall then
take
\begin{equation}\label{eqn:orthorientsigma'}
  \mathcal{B}_{\Sigma'}=\left \{
\begin{pmatrix}
    \mathbf{v} \\
    \mathrm{Jac}_x(\alpha)\mathbf{v} \\
    \mathbf{v} \\
    -\mathrm{Jac}_x(\alpha)\mathbf{v} \
  \end{pmatrix},
\begin{pmatrix}
    \overrightarrow{\mathbf{\xi ^\sharp}}\\
    \mathrm{Jac}_x(\alpha)\overrightarrow{\mathbf{\xi ^\sharp}} \\
    \overrightarrow{\mathbf{\xi ^\sharp}} \\
    -\mathrm{Jac}_x(\alpha)\overrightarrow{\mathbf{\xi ^\sharp}} \
  \end{pmatrix},
\begin{pmatrix}
    \mathbf{\kappa} /\sqrt 2\\
    \mathrm{Jac}_x(\alpha)\mathbf{\kappa}/\sqrt 2 \\
    \mathbf{\kappa} /\sqrt 2\\
    -\mathrm{Jac}_x(\alpha)\mathbf{\kappa}/\sqrt 2 \
  \end{pmatrix},
\begin{pmatrix}
    0 \\
    \alpha _x/\sqrt 2 \\
    0 \\
    -\alpha _x/\sqrt 2 \
  \end{pmatrix}
\right \}.
\end{equation}
Let us set $\varphi _{JX}=\frac12 \left
[\mathrm{Jac}_x(\alpha)^t-\mathrm{Jac}_x(\alpha)\right
]\mathbf{JX}$. We may extend $\mathcal{B}_{\Sigma'}$ to a basis
$\mathcal{B}_{\Gamma}=\mathcal{B}_{\Sigma'}\cup \mathcal{B}_{N}$
of $T_{\zeta}(\Gamma _{\tilde \jmath\times \tilde \jmath})$ by
letting
\begin{equation}\label{eqn:complementsigmabyn}
  \begin{array}{ccl}
  \mathcal{B}_{N}& =&\left \{
\begin{pmatrix}
    0 \\
    \alpha _x /\sqrt 2\\
    0 \\
    \alpha _x /\sqrt 2\
  \end{pmatrix},
\begin{pmatrix}
    \mathbf{v} \\
    \mathrm{Jac}_x(\alpha)^t\mathbf{v} \\
    \mathbf{v} \\
    -\mathrm{Jac}_x(\alpha)^t\mathbf{v} \
  \end{pmatrix},
\begin{pmatrix}
    \overrightarrow{\mathbf{\xi ^\sharp}} \\
    \mathrm{Jac}_x(\alpha)^t\overrightarrow{\mathbf{\xi ^\sharp}} \\
    \overrightarrow{\mathbf{\xi ^\sharp}} \\
    -\mathrm{Jac}_x(\alpha)^t\overrightarrow{\mathbf{\xi ^\sharp}} \
  \end{pmatrix},

\right.\\
&&\left .\begin{pmatrix}
    \mathbf{v} \\
    \mathrm{Jac}_x(\alpha)^t\mathbf{v} \\
    -\mathbf{v} \\
    \mathrm{Jac}_x(\alpha)^t\mathbf{v} \
  \end{pmatrix},
\begin{pmatrix}
    \mathbf{v} \\
    \mathrm{Jac}_x(\alpha)\mathbf{v} \\
    -\mathbf{v} \\
    \mathrm{Jac}_x(\alpha)\mathbf{v} \
  \end{pmatrix},
\begin{pmatrix}
    \overrightarrow{\mathbf{\xi ^\sharp}} \\
    \mathrm{Jac}_x(\alpha)\overrightarrow{\mathbf{\xi ^\sharp}} \\
    -\overrightarrow{\mathbf{\xi ^\sharp}} \\
    \mathrm{Jac}_x(\alpha)\overrightarrow{\mathbf{\xi ^\sharp}} \
  \end{pmatrix},\begin{pmatrix}
    \overrightarrow{\mathbf{\xi ^\sharp}} \\
    \mathrm{Jac}_x(\alpha)^t\overrightarrow{\mathbf{\xi ^\sharp}}\\
    -\overrightarrow{\mathbf{\xi ^\sharp}} \\
    \mathrm{Jac}_x(\alpha)^t\overrightarrow{\mathbf{\xi ^\sharp}} \
  \end{pmatrix},\right.\\
&&\left.
\begin{pmatrix}
    0\\
    \varphi _{\overrightarrow{\mathbf{J\xi ^\sharp}}} \\
    0 \\
    0 \
  \end{pmatrix},\begin{pmatrix}
    0\\
     0\\
    0 \\
    \varphi _{\overrightarrow{\mathbf{J\xi ^\sharp}}} \
  \end{pmatrix},\begin{pmatrix}
    \mathbf{\kappa}/\sqrt 2 \\
    \mathrm{Jac}_x(\alpha)\mathbf{\kappa}/\sqrt 2 \\
    -\mathbf{\kappa} \\
    \mathrm{Jac}_x(\alpha)\mathbf{\kappa}/\sqrt 2\
  \end{pmatrix}
\right \}
  \end{array}
\end{equation}
Performing appropriate column operations, we obtain a basis whose
first $4n+2$ vectors lie in $T_{\zeta}^*(X\times X)$, while the
remaining ones project down to an oriented orthonormal  basis in
$T_{(x,x)}(X'\times X')$. The statement follows from this and
(\ref{eqn:halfformmeta}).

\bigskip

This holds at any point where $r=1$. Taking into account the
appropriate homogeneity, on the upshot we have proved:

\begin{prop}\label{prop:symboljperj}
The symbol of $(\jmath \times \jmath)^*\left (\tilde \Pi
_{X}\right )$ is the section of $\mathrm{Spin}(\Sigma')=\bigwedge
^{1/2}(\Sigma')\otimes \mathcal{S}(N_U)$ given by
\begin{equation}
  \begin{array}{ccl}
    \sigma \left ((\jmath \times \jmath)^*\left (\tilde \Pi _{X}\right
)\right )(x,r\alpha '_{x'},x,-r\alpha '_{x}) &
=&2^{g}r^{g/2}\,\mathrm{vol}_{\Sigma'}^{1/2}\otimes \sigma
_{\frak{g}^M,J}\otimes \sigma _{J_0} \\
     & &((x,r\alpha '_{x},x,-r\alpha' _{x})\in \Sigma ').
  \end{array}
\end{equation}
\end{prop}

\section{The symbol of $\tilde \Pi _{X|X_0}$.}
We shall now examine the symbol of $\tilde \Pi _{X|X_0}$, and
relate it to the symbol of the Szeg\"{o} kernel of $X_0$, $\sigma
(\tilde \Pi _{X_0})\in S^{1/2}(\Sigma _0)$.

\begin{thm}
The symbol of $\tilde \Pi _{X|X_0}$, denoted $\sigma (\tilde \Pi
_{X|X_0})\in S^{(1+g)/2}(\Sigma_0)$, is given by
  \begin{eqnarray*}
    \sigma (\tilde \Pi_{X|X_0}) (\hat{x}, r\alpha _{0\,\hat{x}},\hat{x},-r\alpha _{0\,\hat{x}})
     &=& r^{g/2}\, V_{\mathrm{eff}}(\hat{x})\,\sigma (\tilde \Pi _{X_0})
(\hat{x}, r\alpha _{0\,\hat{x}},\hat{x},-r\alpha
_{0\,\hat{x}})\\
& & (\hat{x}\in X_0,\, r>0).
  \end{eqnarray*}
\label{thm:symbol}
\end{thm}

\noindent \textit{Proof.} We shall now be more sketchy, since the
proof is similar to that of Proposition \ref{prop:symboljperj}. By
homogeneity, it suffices to prove the Theorem for $r=1$.

Let as above $p_{X'}:X'\rightarrow X$ be the projection. The
relevant canonical relation is now $\Gamma _{p_{X'}\times p_{X'}}$
in (\ref{defn:canrelpperp}), and it is naturally diffeomorphic to
the horizontal cotangent bundle of $p_{X'}\times p_{X'}$. The
fibre product $F_{p_{X'}\times p_{X'}}$ of $\Sigma'$ and $\Gamma
_{p_{X'}\times p_{X'}}$ is naturally diffeomorphic to $\Sigma'$.
The projection $p_{F_{p_{X'}\times p_{X'}}}:F_{p_{X'}\times
p_{X'}}\cong \Sigma'\rightarrow \Gamma _{p_{X'}\times p_{X'}}\circ
\Sigma'=\Sigma_0$ is now the quotient map under the $G$-action.

Let $q':T^*(X_0\times X_0)\times T^*(X'\times X')\rightarrow
T^*(X'\times X')$ be the projection onto the second factor. Fix
for now $\zeta '=(x,r\alpha'_x,x,-r\alpha'_x)\in \Sigma'$; the
differential of $q'$ at $(\zeta_0,\zeta')$,
$$q'_{\zeta'}=:d_{(\zeta_0,\zeta')}q: T_{\zeta_0}\left (T^*(X_0\times
X_0)\right )\times T_{\zeta'}\left (T^*(X'\times X')\right )
\rightarrow T_{\zeta'}\left (T^*(X'\times X')\right ),$$ is simply
projection onto the second factor. To ease the exposition, we
shall now assume that the local coordinates on the open
neighbourhood $U\subseteq X'$ of $x$ have been so chosen that the
vector fields $\xi ^\sharp$ are constant (since the action of $G$
along $\Phi ^{-1}(0)$ is free, this may certainly be done).

\begin{lem}\label{lem:localvanishing}
In this system of local coordinates, we have
$$\mathrm{Jac}_y(\alpha')\xi ^\sharp(y)=\mathrm{Jac}_y(\alpha')^t\xi
^\sharp(y)=0,$$ for any $y\in U$, $\xi \in \frak{g}$.\end{lem}

\noindent \textit{Proof.} The two vanishings are equivalent, since
the vector fields $\xi^\sharp$, $\xi \in \frak{g}$, span the
kernel of $d\alpha '$. Since in addition the $\xi ^\sharp$'s are
horizontal on $X'$, the asserted vanishings follow by
differentiating the equality $\alpha '(\xi ^\sharp)=0$.

\bigskip

In these local coordinates,
\begin{equation*}
\begin{array}{ccl}
  q'_{\zeta'}\left (T _{(\zeta_0,\zeta')}(\Gamma _{\jmath\times \jmath})\right ) & =&\left \{
\begin{pmatrix}
  v \\
  \phi \\
 w \\
  \eta
\end{pmatrix}\,:\,v,w ,\phi,\eta \in
\mathbb{R}^{2n-g+1}\text{ such that }\right.\\
&&\left . \phi (\xi ^\sharp (x))=\eta (\xi ^\sharp (x))=0, \,
\forall \,\xi \in \frak{g}\right \}.
\end{array}
\end{equation*}
The symplectic annihilator in $T_{\zeta'}\big (T^*(X)\times
T^*(X)\big )$ is then
\begin{equation*}
\begin{array}{ccl}
  q_{\zeta'}\left (T _{(\zeta',\zeta)}(\Gamma _{\jmath\times \jmath})\right )^{\perp} & =&\left \{
\begin{pmatrix}
  \xi ^\sharp (x)\\
 0 \\
 \eta ^\sharp (x) \\
0
\end{pmatrix}\,:\,\xi, \eta \in \frak{g}\right
\}.
\end{array}
\end{equation*}
We shall now denote by $V_0$, $V_1$ and $V$ the analogues in the
present setting of the vector bundles $U_0$, $U_1$ and $U$
introduced in section 3. In view of Lemma \ref{lem:localvanishing}
we now have
\begin{equation*}
\begin{array}{ccl}
  V_{0\zeta'} & =&\left \{
\begin{pmatrix}
  \xi ^\sharp (x)\\
 0 \\
 \xi ^\sharp (x) \\
0
\end{pmatrix}\,:\,\xi \in \frak{g}\right
\}\text{ and }\,\,\,\,\, V_{1\zeta}=\left \{\begin{pmatrix}
  \xi ^\sharp (x)\\
 0 \\
 \eta ^\sharp (x) \\
0
\end{pmatrix}\,:\,\xi, \eta \in \frak{g}\right
\}.
\end{array}
\end{equation*} Hence, in the terminology of Corollary \ref{cor:imageEpmsurjmorph}
iii) and recalling Corollary \ref{cor:normalbundleU}, we obtain
for the quotient bundle $V=V_1/V_0\subseteq N_{\Sigma'}\cong N_U$:
$$V_{\zeta'}=\left \{\begin{pmatrix}
  \xi ^\sharp (x)\\
 0 \\
 -\xi ^\sharp (x) \\
0
\end{pmatrix}\,:\,\xi\in \frak{g}\right
\}\cong \frak{g}^M_r(p)\subseteq N_U.
$$
Again by Corollary \ref{cor:imageEpmsurjmorph}, the symplectic
normal bundle $N_{V}$ of $V$ in $N_{\Sigma'}$ is thus
$$N_{V}
\cong (q^*H(M'/M_0),2r\Omega
 _{H(M/M')})\oplus (q^*H(M'/M_0),-2r\Omega _{H(M/M')}).$$
This is clearly simplectically isomorphic to the normal bundle of
$\Sigma _0$ in $T^*(X_0\times X_0)$ (cfr Proposition 6.4 of
\cite{bg}). By Corollary \ref{cor:imageEpmsurjmorph} v), the image
of $\sigma _{\frak{g}^M,J}\otimes \sigma _{J_0}$ under the vector
bundle morphism $\Phi _{V}:\mathcal{S}(N_{\Sigma'})\rightarrow
\bigwedge ^{-1/2}(V)\otimes \mathcal{S}(N_{V})$ is
$$\Phi _{V}(\sigma _{\frak{g}^M,J}\otimes \sigma
_{J_0})=\mathrm{vol}_{V}^{-1/2}\otimes \sigma _{J_0},$$ where
$\mathrm{vol}_{V}^{-1/2}$ is the $-\frac 12$-form on the oriented
isotropic subbundle $V\subseteq T(T^*(X'\times X'))$ taking value
one on oriented orthonormal basis of $V$.

The following Lemma is proved with an argument similar to that
used for Lemma \ref{cor:normalbundleU}, the main difference being
we need to include an oriented orthonormal basis for
$V_{0,\zeta'}$.

\begin{lem}\label{lem:normalbundleU}
Suppose $x\in X'$ and let $\zeta '=(x,\alpha' _x,x,-\alpha '_x)$
(thus, $r=1$). Then the image at $\zeta'$ of
$\mathrm{vol}_{\Gamma_{p_{X'} \times p_{X'}},\zeta'} ^{1/2}\otimes
\mathrm{vol}_{\Sigma',\zeta '}^{1/2}\otimes
\mathrm{vol}_{V,\zeta'} ^{-1/2}$ in $\bigwedge
^{1/2}(\Sigma_0)_{\zeta _0}\otimes \det (V_{0,\zeta'}^{ *})$ under
the line bundle isomorphism (\ref{eqn:isomolnbdl}) is
$2^{-g}\,\mathrm{vol}_{\Sigma_0,\zeta_0}^{1/2}\otimes
\mathrm{vol}_{V_{0,\zeta'}}$, where $\mathrm{vol}_{V_{0}}$ is the
volume form on $V_{0}$ taking value one on oriented orthonormal
frames.
\end{lem}

Clearly, $\mathrm{vol}_{V_{0,\zeta'}}$ is the volume form on the
vertical tangent bundle of the principal $G$-bundle $\Sigma
'\rightarrow \Sigma _0$ associated to the orientation and the
restricted metric. The statement of Theorem \ref{thm:symbol} then
follows by fibrewise integration and homogenization.

\section{The Asymptotic Expansion}

The function $\sum _{j=1}^{d_k}s_j^{(G,k)}(x')\otimes
\overline{s_j^{(G,k)}(y')}$ appearing in (\ref{eqn:localformpf})
is obviously well-defined on $X_0$. We shall now argue that it
admits an asymptotic expansion as in the statement of Theorem
\ref{thm:mainexpgequiv}.

To this end, let us introduce the following auxiliary
Fourier-Hermite distribution on $X_0\times X_0$:
\begin{equation}\label{eqn:auxiliary}
\begin{array}{ccl}
 \tilde P_{X|X_0}&=:&(V_{\mathrm{eff}}\boxtimes V_{\mathrm{eff}})^{-1}\,\tilde \Pi _{X|X_0}\\
 &=&\left (\sum _{k=0}^{+\infty}\sum _{j=1}^{d_k}s_j^{(G,k)}\boxtimes \overline{s_j^{(G,k)}}\right )\mathrm{vol}^{1/2}_{X_0}\boxtimes
\mathrm{vol}^{1/2}_{X_0}.\\
\end{array}
\end{equation}
Thus $\tilde P_{X|X_0}\in J^{(1+g)/2}(X_0\times X_0,\Sigma _0)$,
and the associated operator
$P_{X|X_0}:\mathcal{D}'_{1/2}(X_0)\rightarrow
\mathcal{D}'_{1/2}(X_0)$ is an $S^1$-invariant elliptic operator
satisfying $\Pi _{X_0}\circ P_{X|X_0}=P_{X|X_0}$. One can now
follow the arguments in \cite{stz}, Lemma 4.2 and 4.3 and \cite{g}
to conclude that $P_{X|X_0}$, and thus $\Pi _{X|X_0}$, is an
elliptic Toeplitz operator possessing a semiclassical symbol, and
that this implies the asserted asymptotic expansion.

\end{document}